\documentstyle[12pt]{article}
\input amssym.def
\begin {document}
\def \Z{\Bbb Z}
\def \C{\Bbb C}

\def \Q{\Bbb Q}

\def \wt{{\rm wt}\;}
\def \fg{\frak g}

\def \Res{{\rm Res}}
\def \End{{\rm End}\;}

\def \Irr {{\rm Irr}}

\def \Hom{{\rm Hom}}
\def \mod{{\rm mod}}

\def \<{\langle} 
\def \>{\rangle} 

\def \a{\alpha }
\def \e{\epsilon }

\def \g{\gamma}
\def \b{\beta }

\def \be{\begin{equation}\label}
\def \ee{\end{equation}}
\def \bl{\begin{lem}\label}
\def \el{\end{lem}}
\def \bt{\begin{thm}\label}
\def \et{\end{thm}}
\def \bp{\begin{prop}\label}
\def \ep{\end{prop}}
\def \br{\begin{rem}\label}
\def \er{\end{rem}}
\def \bc{\begin{coro}\label}
\def \ec{\end{coro}}
\def \bd{\begin{de}\label}
\def \ed{\end{de}}
\def \pf{{\bf Proof. }}

\newtheorem{thm}{Theorem}[section]
\newtheorem{prop}[thm]{Proposition}
\newtheorem{coro}[thm]{Corollary}

\newtheorem{lem}[thm]{Lemma}
\newtheorem{rem}[thm]{Remark}
\newtheorem{de}[thm]{Definition}

\makeatletter
\@addtoreset{equation}{section}

\makeatother
\makeatletter

\baselineskip=16pt
\begin{center}{\Large \bf Certain
extensions of vertex operator algebras of affine type}

\vspace{0.5cm}
Haisheng Li\footnote{Partially supported by NSF grant
DMS-9970496}\\ Department of Mathematical Sciences, Rutgers University,
Camden, NJ 08102
\end{center}

{\bf Abstract} We generalize Feigin and Miwa's construction
of extended vertex operator (super)algebras $A_{k}(sl(2))$ 
for other types of simple Lie algebras. For all the constructed extended 
vertex operator (super)algebras, irreducible modules are 
classified, complete reducibility of every module is proved
and fusion rules are determined modulo the fusion rules for
vertex operator algebras of affine type.

\section{Introduction}

In the development of vertex operator algebra theory,
one of the most important problems is to construct 
new {\em solvable} vertex operator (super)algebras 
in the sense that irreducible modules and fusion rules can be 
completely determined and that intertwining operators 
can be explicitly constructed. To a certain extent, 
such algebras give rise to solvable physical models.
One of many ways to get such vertex operator (super)algebras 
is to consider certain extensions of some well known algebras.
For example ([MS], [Li5]), the vertex operator algebra 
$L(k,0)$ associated to the affine Lie algebra $\hat{s}l_{2}$ 
with a positive {\em even} integral level $k$ can be extended 
to a vertex operator (super)algebra $L(k,0)+L(k,k)$. 
When $k$ is odd, $L(k,0)+L(k,k)$ does not 
have an extended vertex operator (super)algebra 
structure because of the failure of the locality. 
(For a certain class of vertex operator algebras,
e.g., vertex operator algebras associated to 
positive-definite even lattices, as proved in [DL], 
the sum of a copy of each irreducible module
does have a nice structure, called an abelian intertwining algebra.
See [DL], [FFR], [M] and [Hua2] for notions of various
generalized structures.)

In [FM], Feigin and Miwa constructed a family of extended vertex operator 
(super)algebras $A_{k}$ from the vertex operator algebras $L(k,0)$ 
associated to the affine Lie algebra $\hat{s}l_{2}$ with an {\em arbitrary}
positive integral level $k$, and they classified all irreducible modules
and determined all fusion rules. 
In addition they obtained very interesting results
on the monomial basis for irreducible modules.
This paper was mainly motivated by [FM] and [DLM2].
As the main results of this paper we generalize their results
except the monomial basis result to affine Lie algebras of other 
types by using a different approach.

The algebras $A_{k}$ were defined in [FM] 
by a set of mutually local vertex operators (or fields).
On the other hand, in terms of vertex operator algebra language,
$A_{k}$ are extensions (by an infinite sum of irreducible modules) of 
vertex operator algebra
$L(k,0)\otimes M(1,0)$, where $M(1,0)$ is the vertex operator algebra
associated to an infinite-dimensional Heisenberg Lie algebra of rank one, 
or a single free bosonic field. The essential building block of $A_{k}$ is 
the irreducible $L(k,0)\otimes M(1,0)$-module
$L(k,k)\otimes M(1,\a)$ for some $\a\in \C$.

The $L(k,0)$-module $L(k,k)$ has been known to be a {\em simple current} 
([FG], [GW], [SY]) in the sense that
the left multiplication of the equivalence class
$[L(k,k)]$ in the Verlinde algebra gives rise to a permutation on
the standard basis.
It was known to physicists (cf. [MS], [SY]) that a simple current with 
integer weights can be included to generate an extended 
vertex operator algebra.
In [Li5], as an exercise by using an explicit construction of
simple currents given in [Li4] we studied the 
extension of a certain vertex operator algebra by a self-dual
(or order 2) simple current where $L(k,0)+L(k,k)$
is a special case. For such extended algebras,
all their irreducible modules were classified and the complete
reducibility of every module was proved.
A little bit latter, the results of [Li5] were 
greatly extended in [DLM2].

The construction of simple currents given in [Li4] 
was based on a result of [Li2].
Let $V$ be a vertex operator algebra and
let $h$ be a weight one primary vector
in $V$ such that the component operators $h(m)$ of the vertex operator 
$Y(h,z)$ satisfy the Heisenberg
algebra relation and such that $h(0)$ is semisimple on $V$ with only
rational eigenvalues. Clearly, $\sigma_{h}:=e^{2\pi ih(0)}$ is an automorphism 
of $V$. Define
$$\Delta(h, z)
=z^{h(0)}\exp \left(\sum_{n=1}^{\infty}{h(n)\over -n}
(-z)^{-n}\right),$$
an element of $(\End V)\{z\}$.
It was proved in [Li2] that for any $V$-module 
$W$,
$$(W^{(h)},Y_{h}(\cdot,z)):= (W,Y(\Delta(h,z)\cdot,z))$$ 
is a $\sigma_{h}$-twisted $V$-module.
In particular, this gives an (untwisted) $V$-module
if $h(0)$ acting on $V$ has only integral eigenvalues.
It was furthermore proved in [Li4] that if $I(\cdot,z)$ is an intertwining 
operator of type ${W_{3}\choose W_{1}W_{2}}$ (in the sense of [FHL]),
then $I(\Delta(h,z)\cdot,z)$ is an intertwining 
operator of type ${W_{3}^{(h)}\choose W_{1}W_{2}^{(h)}}$.
By using this result and the invertiability of $\Delta(h,z)$,
it was proved that $V^{(h)}$ is a simple current. 
As a matter of fact, for certain well known vertex operator algebras 
almost all simple currents can be constructed in this way.
For example, when $V$ is the vertex operator algebra $V_{L}$ associated with
a positive-definite even lattice $L$, it was proved that
all irreducible $V_{L}$-modules can be constructed this way.
In this case, this construction is intimately related to the construction
of twisted modules by using shifted vertex operators in [Le].
When $V=L(\ell,0)$ associated 
to an affine Lie algebra $\hat{\fg}$ with $\fg\ne E_8$ and with a positive 
integral level $\ell$, all the simple currents
can also be constructed this way. 
The merit of this construction of a simple current
is on the canonicalness of the vector space and the vertex operator map
(in terms of the algebra $V$ and the element $h$). 
With this construction, certain intertwining operators can be constructed 
canonically. The canonical construction of
intertwining operators of certain types
is the basis of [Li5] and [DLM2].

The essential results of [DLM2] can be described as follows:
Let $V$ be a vertex operator algebra and
$H$ be a subspace of $V_{(1)}$ such that the components
$h(n)$ of vertex operators $Y(h,z)$ for $h\in H,\; n\in \Z$
satisfy the Heisenberg algebra relation. Let $L$ be a subgroup of $H$ 
such that for every $\a\in L$, $\a(0)$ acts semisimply on $V$ with only
integral eigenvalues. Consider the space $V[L]=\C[L]\otimes V$,
where for $\a\in L$, $e^{\a}\otimes V$ is identified with $V^{(\a)}$ 
equipped with the $V$-module structure $Y_{\a}$.
Extend the $V$-module structure on $V[L]$ 
in a certain canonical way
to a vertex operator map $Y$ on $V[L]$. Then it was proved that $V[L]$ 
equipped with the defined vertex operator map $Y$ is a 
generalized vertex algebra in the sense of [DL]. 
It was proved that $V[L]$ is a vertex operator (super)algebra
when $L$ satisfies certain conditions.
When $V=M_{\bf h}(1,0)$ with ${\bf h}=\C\otimes _{\Z}L$, where $L$ is an
integral lattice, we have $V[L]=V_{L}$ ([FLM], [DL]). 
Then we may view $V[L]$ as a 
generalization of $V_{L}$.
In [DLM2], we were mainly interested in the case $V=L(k,0)$
associated to an affine algebra $\hat{\fg}$ with a positive integral level $k$.
Because $L(k,0)$ has only finitely many irreducible modules up to 
equivalence, each irreducible $V$-module $V^{(\a)}$ in $V[L]$ is not
multiplicity-free. Having noticed this we proved that 
$V[L]$ has a quotient algebra $\overline{V[L]}$ such that
every irreducible $V$-module in $\overline{V[L]}$ is 
multiplicity-free and that $V[L]$ and $\overline{V[L]}$ contain
the same number of non-isomorphic irreducible $V$-modules.
An irreducible $V$-module $W$ was also extended
to $W[L]$ on which $V[L]$ acts and it was proved that $W[L]$ is in general
a twisted $V[L]$-module with respect to a certain 
automorphism of $V[L]$. With this result, under the assumption
of complete reducibility of $V$-modules, all irreducible
$V[L]$-modules were classified and a complete reducibility
theorem for $V[L]$-modules was proved.

In this paper, to generalize Feigin and Miwa's construction
we apply the results of [DLM2] by taking
$V=L(k,0)\otimes M_{{\bf h}'}(1,0)$ and by choosing ${\bf h}'$ and
$L$ appropriately, depending on the type of $\fg$. In this case, 
each irreducible $V$-module in $V[L]$ is multiplicity-free and
$V[L]$ is a simple algebra.
On the other hand, since the category of $V$-modules is not semisimple,
the results of [DLM2] for the complete reducibility of $V[L]$-modules
do not apply to this case directly.
The complete reducibility of $V[L]$-modules is proved here.
We also naturally  extend
intertwining operators for modules in the category of 
$V$-modules to intertwining operators for modules 
in the category of $V[L]$-modules. Using this result 
we are able to derive a formula of fusion rules for 
$V[L]$-modules in terms of fusion rules for $V$-modules.

This paper is organized as follows: In Section 2,
we recall the classification of irreducible modules for certain 
vertex operator algebras and recall a construction
of simple currents. Most part of this section is preliminary and
the only new result is about
the fusion rules for simple currents for $L(k,0)$.
In Section 3 we recall and refine some of the results of [DLM2] on 
the extended vertex (super)algebra $V[L]$. We furthermore study
the multiplicity-free case.
In Section 4, we apply the results of Section 3 to
construct extended vertex 
operator (super)algebras associated to an affine algebra $\hat{\fg}$.

\section{Vertex operator algebras and simple currents}

The extended vertex operator (super) algebras we shall construct are
based on vertex operator algebras 
$L_{\fg}(\ell,0)$, $M_{\bf h}(1,0)$, $V_{L}$ and their
representations. For this reason, in this section 
we shall recall the relevant information about these algebras and
we also recall from [Li4] a construction of 
simple currents for a certain type of vertex operator algebras, 
including $L_{\fg}(\ell,0)$, $M_{\bf h}(1,0)$, $V_{L}$.
For the vertex operator algebra $L_{\fg}(\ell,0)$ associated to 
an affine algebra $\hat{\fg}$ (not type $E_{8}^{(1)}$)
with a positive integral level $\ell$, we prove that the equivalence 
classes of the simple currents form an abelian group isomorphic to
$P^{\vee}/Q^{\vee}$, where $P^{\vee}$ and $Q^{\vee}$ are the co-weight 
and co-root lattices of $\fg$.

\subsection{Vertex operator algebras $L_{\fg}(\ell,0)$, $M_{\bf h}(1,0)$
and $V_{L}$}

We shall use standard definitions and notations 
as given in [FHL] and [FLM] and we also use [K] and [H]
as our references for (Kac-Moody) Lie algebras.
Following [DL] we use the term ``vertex (super)algebra''
for an object that satisfies all the axioms defining the notion of
vertex operator (super)algebra except the two grading restrictions.


Let ${\fg}$ be a finite-dimensional simple Lie algebra,
${\bf h}$ a Cartan subalgebra, and $\<\cdot,\cdot\>$
the normalized killing form such that the square length of a long root is 2.
Let 
\begin{eqnarray}
\hat{\fg}=\fg \otimes \C[t,t^{-1}]\oplus \C c
\end{eqnarray}
be the affine Lie algebra.

Let $\ell$ be a complex number such that
$\ell\ne -h^{\vee}$, where $h^{\vee}$ is the dual Coxeter number of $\fg$.
Let $\C_{\ell}$ be the one-dimensional $(\fg \otimes \C[t]+\C c)$-module
on which $c$ acts as scalar $\ell$ and $\fg\otimes \C[t]$ 
acts as zero. Form the generalized Verma $\hat{\fg}$-module
\begin{eqnarray}
M_{\fg}(\ell,0)=U(\hat{\fg})\otimes_{U(\fg\otimes \C[t]+\C c)}\C_{\ell}.
\end{eqnarray}
It was well known (cf. [FF], [FZ], [Li1], [MP]) 
that $M_{\fg}(\ell,0)$ has a natural
vertex operator algebra structure. Furthermore, the category
of {\em weak} $M_{\fg}(\ell,0)$-modules in the sense that all the axioms 
defining the notion of module except those involving grading hold is 
canonically equivalent to the category of {\em restricted} (cf. [K])
$\hat{\fg}$-modules
of level $\ell$ in the sense that for every vector $w$ of the module,
$(\fg\otimes t^{n}\C[t])w=0$ for $n$ sufficiently large.

\br{rvoageneral}
{\em More generally, let $W$ be an arbitrary vector space. 
An element $a(z)$ of $(\End W)[[z,z^{-1}]]$ is called a
{\em vertex operator} if $a(z)w\in W((z))$ for every $w\in W$.
Two vertex operators $a(z)$ and $b(z)$ are said to be 
{\em mutually local} if there exists a nonnegative integer $N$ such that
\begin{eqnarray}\label{elocal}
(z_{1}-z_{2})^{N}[a(z_{1}),b(z_{2})]=0
\end{eqnarray}
(cf. [DL], (1.4)).
It was proved in [Li1] (Corollary 3.2.11) that
any set of mutually local vertex operators on $W$ automatically 
generates a vertex algebra, which is a canonical vector subspace of 
$(\End W)[[z,z^{-1}]]$, and that $W$ is a natural module for 
this vertex algebra.}
\er

For $\lambda\in {\bf h}^{*}$, denote by $L_{\fg}(\ell,\lambda)$ 
the irreducible highest weight
$\hat{\fg}$-module of level $\ell$ with highest weight $\lambda$.
Each $L_{\fg}(\ell,\lambda)$ is an irreducible $M_{\fg}(\ell,0)$-module
possibly with infinite-dimensional homogeneous subspaces.

Denote by $\theta$ the highest long root of $\fg$.
For a positive integer $\ell$, set
\begin{eqnarray}
P_{\ell}=\{\lambda\in P_{+}\;|\; \<\lambda, \theta\>\le \ell\},
\end{eqnarray}
where $P_{+}$ is the set of dominant integral weights of $\fg$.
Then $L(\ell,\lambda)$ is an integrable $\hat{\fg}$-module if and only 
if $\lambda\in P_{\ell}$ [K].
The following result was known (cf. [DL, Proposition 13.17], 
[FZ, Theorem 3.1.3], [Li1, Propositions 5.2.4 and 5.2.5], 
[MP, Theorems 5.9, 5.14 and 5.15]):

\bp{pFZDL}
Let $\ell$ be a positive integer. Then 
(1) The set of irreducible $L_{\fg}(\ell,0)$-modules is exactly
the set of irreducible highest 
weight integrable (or standard) $\hat{\fg}$-modules of level $\ell$.
(2) Every $L_{\fg}(\ell,0)$-module is completely reducible.
\ep

The following stronger result was obtained in [DLM1] (Theorem 3.7):

\bp{pDLM}
Let $\ell$ be a positive integer. Then every weak 
$L_{\fg}(\ell,0)$-module is a direct sum of irreducible highest 
weight integrable (or standard) $\hat{\fg}$-modules of level $\ell$. 
In particular, every irreducible weak $L_{\fg}(\ell,0)$-module is 
an (ordinary) $L_{\fg}(\ell,0)$-module.
\ep

A vertex operator algebra with the property that every weak module is a 
direct sum of irreducible (ordinary) modules is said to be {\em regular} 
[DLM1]. 

Let ${\bf h}$ be a finite-dimensional vector space equipped with
a nondegenerate symmetric bilinear form $\<\cdot,\cdot\>$, i..e,
a finite-dimensional abelian Lie algebra equipped with
a nondegenerate symmetric invariant bilinear form.
Let $\hat{\bf h}={\bf h}\otimes \C[t,t^{-1}]+\C c$ be the affine
Lie algebra. We have
\begin{eqnarray}
\hat{\bf h}=\hat{\bf h}_{\Z}\oplus {\bf h},
\end{eqnarray}
where $\hat{\bf h}_{\Z}=\sum_{n\ne 0}{\bf h}\otimes t^{n}+\C c$ is a 
Heisenberg algebra and ${\bf h}$ is central.
Similar to the construction of $M_{\fg}(\ell,0)$ we construct 
a space $M_{\bf h}(1,0)$, and
just like $M_{\fg}(\ell,0)$, $M_{\bf h}(1,0)$ is a 
vertex operator algebra. For any $\a\in {\bf h}^{*}$ $(={\bf h})$,
let $\C e^{\a}$ be a one-dimensional $({\bf h}\otimes \C[t]+\C c)$-module on
which ${\bf h}\otimes t\C[t]$ acts as zero, $c$ acts as $1$ and
$h=h(0)$ acts as scalar $\<\a,h\>$ for $h\in {\bf h}$.
Form the induced module
\begin{eqnarray}
M_{\bf h}(1,\a)=U(\hat{\bf h})\otimes_{U({\bf h}\otimes \C[t]+\C c)}\C e^{\a}
\simeq M_{\bf h}(1,0)\otimes \C e^{\a}\;\;\;\mbox{(linearly)}.
\end{eqnarray}
As in the case of $M_{\fg}(\ell,0)$, 
$M_{\bf h}(1,\a)$ is an irreducible $M_{\bf h}(1,0)$-module.
Furthermore, from [FLM] the lowest $L(0)$-weight of $M_{\bf h}(1,\a)$ is 
\begin{eqnarray}\label{elowestweightlattice}
\Delta_{\a}={1\over 2}\<\a,\a\>.
\end{eqnarray}
On the other hand, clearly every irreducible 
$M_{\bf h}(1,0)$-module is isomorphic to $M_{\bf h}(1,\a)$
for some $\a$. It follows from the complete reducibility
of certain $\hat{\bf h}_{\Z}$-modules ([LW], [K]) that
an $M_{\bf h}(1,0)$-module on which ${\bf h}$ semisimply acts
is completely reducible.
In general, an $M_{\bf h}(1,0)$-module may not be completely reducible.

Let $P$ be a rational lattice of finite rank with the $\Z$-bilinear 
form $\<\cdot,\cdot\>$. 
Set
\begin{eqnarray}
{\bf h}=\C\otimes_{\Z}P,
\end{eqnarray}
and extend $\<\cdot,\cdot\>$ to a $\C$-bilinear form on ${\bf h}$.

Denote by $\C[P]$ the group algebra. Set
\begin{eqnarray}
V_{P}=\C[P]\otimes M_{\bf h}(1,0),
\end{eqnarray}
equipped with the standard $M_{\bf h}(1,0)$-module 
(or $\hat{\bf h}_{\Z}$-module) structure. That is, $V_{P}$ is
a direct sum of irreducible 
$M_{\bf h}(1,0)$-modules 
$M_{\bf h}(1,\a)\simeq \C e^{\a}\otimes M_{\bf h}(1,0)$ 
for $\a\in P$.

It was proved in [FLM] (cf. [B]) that when $P=L$ is even and 
positive-definite, $V_{L}$ has a natural
simple vertex operator algebra structure which extends
the $M_{\bf h}(1,0)$-module structure. (It follows from 
Proposition \ref{puniqueness} that
such a vertex operator algebra structure is unique up to equivalence.)
Furthermore, let $L^{o}$ be the dual lattice of $L.$ Then $V_{L^{o}}$ is a
natural $V_{L}$-module and 
$V_{\beta+L}$ is an irreducible $V_{L}$-module for $\beta\in L^{o}.$

\bp{p} (1) Let $\{\beta_{1},\dots,\beta_{r}\}$ be
a complete set of representatives of cosets of $L$ in $L^{o}$.
Then $\{V_{\beta_{1}+L},\dots, V_{\beta_{r}+L}\}$
is a complete set of representatives of equivalent classes of irreducible
$V_{L}$-modules.
(2) Every $V_{L}$-module is completely reducible.
(3) $V_{L}$ is regular, i.e., every weak $V_{L}$-module is a direct sum of 
irreducible (ordinary) $V_{L}$-modules.
\ep

The assertion (1) was proved in [FLM] and [D1], (2) was proved in [Guo] and
(3) was proved in [DLM1].

\br{rDLgva}
{\em For a general rational lattice $P$, $V_{P}$ is not 
a vertex (operator) algebra
because of the involvement of non-local vertex operator.
A notion of generalized vertex algebra was introduced in [DL] 
with a generalized Jacobi identity as one of the main axioms
and it was proved in [DL, Theorem 5.1 and Remark 9.11] 
that $V_{P}$ is a generalized vertex algebra.}
\er

\subsection{Simple currents and a construction}
We first recall from [FHL] the definition of an intertwining operator.

\bd{dintertwiningoperator}
{\em Let $V$ be a vertex operator algebra and let $W_{1},W_{2}$ and $W_{3}$ 
be $V$-modules.
An {\em intertwining operator } of type ${W_{3}\choose W_{1}W_{1}}$ 
is a linear map $I$ from $W_{1}$ to $(\Hom (W_{2},W_{3}))\{z\}$, 
(where for a vector space $U$, $U\{z\}$ is defined to be the vector
space of $U$-valued formal series in $z$ with arbitrary complex 
powers of $z$), satisfying the following properties:
for $w_{1}\in W_{1},\; w_{2}\in W_{2}$, 
\begin{eqnarray}
I(w_{1},z)w_{2}\in z^{\gamma_{1}}W_{3}[[z]]+\cdots +z^{\gamma_{n}}W_{3}[[z]]
\end{eqnarray}
for some (finitely many) complex numbers $\gamma_{1},\dots,\gamma_{n}$,
and for $v\in V,\; w_{1}\in W_{1}$,
\begin{eqnarray}
[L(-1),I(w_{1},z)]={d\over dz}I(w_{1},z),\label{ederivative}
\end{eqnarray}
\begin{eqnarray}\label{ejacobiinter}
& &z^{-1}_0\delta\left(\frac{z_1-z_2}{z_0}\right)
Y(v,z_1)I(w_{1},z_2)-z^{-1}_0\delta\left(\frac{z_2-z_1}{-z_0}\right)
I(w_{1},z_2)Y(v,z_1)\nonumber\\
&=&z_2^{-1}\delta\left(\frac{z_1-z_0}{z_2}\right)
I(Y(v,z_0)w_{1},z_2).
\end{eqnarray}}
\ed

All intertwining operators of type 
${W_{3}\choose W_{1}W_{2}}$ form a vector space denoted by
$I_{W_{1}W_{2}}^{W_{3}}$. The dimension of $I_{W_{1}W_{2}}^{W_{3}}$
is called the {\em fusion rule}, 
denoted by $N_{W_{1},W_{2}}^{W_{3}}$.
Clearly, the fusion rule only depends on the equivalence class 
of each $W_{i}$.

The following are among the immediate consequences of the Jacobi identity
(\ref{ejacobiinter}):
\begin{eqnarray}\label{ecomm}
[v_{n},I(w_{1},z_{2})]
=\sum_{i\ge 0}{n\choose i}z_{2}^{n}I(v_{i}w_{1},z_2)
\end{eqnarray}
(the {\em commutator formula}) and
\begin{eqnarray}\label{eiterate}
& &I(v_{n}w_{1},z_{2})\nonumber\\
&=&\Res_{z_{1}}\left((z_1-z_2)^{n}
Y(v,z_1)I(w_{1},z_2)-(-z_2+z_1)^{n}I(w_{1},z_2)Y(v,z_1)\right)
\end{eqnarray}
(the {\em iterate formula}) for $n\in \Z$.
Conversely, the commutator and iterate formulas imply the Jacobi identity.

\br{rintertwiningoperator}
{\em If $V=L_{\fg}(\ell,0)$ and $W_{1}=L_{\fg}(\ell,\lambda)$, 
for $w_{1}\in L(\lambda)$
(the lowest  $L(0)$-weight subspace of $L_{\fg}(\ell,\lambda)$), 
the commutator formula (\ref{ecomm}) gives
\begin{eqnarray}\label{ecommaffine}
[a_{n},I(w_{1},z_{2})]=z_{2}^{n}I(aw_{1},z_{2})
\end{eqnarray}
for $a\in \fg\subset L_{\fg}(\ell,0)$ and for $n\in \Z$.
In many literatures such as [TK] (and [FM]), an intertwining operator 
was defined on $L(\lambda)$ with the properties (\ref{ederivative})
and (\ref{ecommaffine}) as the defining axioms. 
However, it can be proved (cf. [TK], [Li3, Li6]) that 
the two definitions give rise to the same fusion rules.}
\er

Let $V$ be a vertex operator algebra. We denote by $\Irr (V)$
the set of equivalence classes of irreducible modules and
for an irreducible $V$-module $W$, denote
by $[W]$ the equivalence class.
$V$ is said to be {\em quasi-rational} [MS]
if all fusion rules associated with irreducible modules
are finite and if for any $[W_{1}], [W_{2}]\in \Irr(V)$,
$N_{[W_{1}],[W_{2}]}^{[W_{3}]}=0$ for all but finitely many $[W_{3}]\in \Irr(V)$.

For a quasi-rational vertex operator algebra $V$, 
the {\em Verlinde algebra} ${\cal{A}}(V)$ is defined
to be an algebra (over $\C$) with $\Irr(V)$ as a 
basis and with the fusion rules as the structural constants, i.e.,
\begin{eqnarray}
[W_{i}]\cdot [W_{j}]
=\sum_{[W_{k}]\in \Irr(V)}N_{[W_{i}],[W_{j}]}^{[W_{k}]}[W_{k}].
\end{eqnarray}
When $V$ is simple, it is easy to show that $[V]$ is the unit.
It follows immediately from [FHL, HL2]
that ${\cal{A}}(V)$ is commutative.
Under certain conditions, Huang [Hua1]
established the associativity,
but for a general $V$
the associativity is still an unsolved problem.

Let $V=L_{sl(2)}(k,0)$ with a positive integer $k$. 
It is well known
([GW], [TK], [FZ]) that the Verlinde algebra has the following relations:
\begin{eqnarray}
[L(k,i)]\cdot [L(k,j)]
=\sum_{r=\mbox{max}(i-j,j-i)}^{\mbox{min}(i+j,2k-i-j)} [L(k,r)].
\end{eqnarray}

The following definition is due to Schellekens and 
Yankielowicz [SY].

\begin{de}\label{d2.10}
{\em 
An irreducible $V$-module $W$ is called a {\em simple current} 
if the associated matrix of the left multiplication of $[W]$
with respect to the standard basis of the Verlinde algebra
is a permutation. The order of the associated matrix as a group element
is called the {\em order} of $W$.}
\end{de}


We now recall a construction of simple currents from [Li4].
In the following, one may think of $V$ as
one of, or more general, any tensor product of the following vertex operator algebras:
$$L_{\fg}(\ell,0),\;\;\; M_{\bf h}(1,0),\;\;\; V_{L},\;\;\;
L_{\fg}(\ell,0)\otimes M_{\bf h}(1,0),\;\;\;
L_{\fg}(\ell,0)\otimes V_{L}.$$

Let $\alpha\in V_{(1)}$ satisfying the following conditions:
\begin{eqnarray}\label{h2.11}
L(n)\alpha=\delta_{n,0}\alpha,\;\;\;
\alpha(n)\alpha=\delta_{n,1}\gamma {\bf 1}\;\;\;\mbox{for }n\in {\Z}_{+},
\end{eqnarray}
where $Y(\a,z)=\sum_{n\in \Z}\a(n)z^{-n-1}$, i.e., $\a(n)=\a_{n}$, and 
$\gamma$ is a fixed complex number. 
Notice that condition (\ref{h2.11}) implies 
that $\a$ is a primary vector and that operators $\a(n)$ satisfy
the Heisenberg algebra relation
\begin{eqnarray}\label{eheisenberg1}
[\a(m),\a(n)]=m\gamma \delta_{m+n,0}\;\;\;\mbox{ for }m,n\in \Z.
\end{eqnarray}
Furthermore, assume that $\a(0)$ acts semisimply on $V$.
It is clear that $e^{2\pi i\a(0)}$ is an automorphism of $V$ and that
$e^{2\pi i\a(0)}=1$ if and only if $\alpha(0)$ has only integral 
eigenvalues on $V$. If each $\a(0)$ has only
rational eigenvalues on $V$ and if $V$ is finitely generated, then 
$e^{2\pi i\a(0)}$ is of finite order.
Define 
\begin{eqnarray}\label{edeltadefinition}
\Delta(\a,z)=z^{\a(0)}\exp\left(\sum_{k=1}^{\infty}\frac{\a(k)}{-k}
(-z)^{-k}\right).
\end{eqnarray}
This is a well defined element of $(\End W)\{z\}$ for any weak $V$-module
$W$ on which $\a(0)$ semisimply acts.
The following result is a special case of Proposition 5.4 of [Li2].

\begin{prop}\label{p2.14} Let $\alpha, \Delta(\alpha,z)$ be given
as before. Assume that $\alpha(0)$ has only integral eigenvalues on $V$.
Let $W$ be any  (irreducible) weak $V$-module. Set
\begin{eqnarray}
(W^{(\a)},Y_{\a}(\cdot,z))=(W, Y(\Delta(\alpha,z)\cdot,z)).
\end{eqnarray}
Then $(W^{(\a)}, Y_{\a})$ carries the structure of
an (irreducible) weak $V$-module.
\end{prop}

As a convention, by $V$-module $W^{(\a)}$ 
we mean the $V$-module $(W^{(\a)},Y_{\a})$.

Recall the following result from [Li4] (Proposition 2.5):

\bp{pli4} Let $\alpha, \Delta(\alpha,z)$ be as in Proposition \ref{p2.14}.
Let $W_{1}$ and $W_{2}$ be weak $V$-modules and
$f$ be a $V$-homomorphism from $W_{1}$ to $W_{2}$. Then
$f$ is also a $V$-homomorphism from $W_{1}^{(\a)}$ to $W_{2}^{(\a)}$.
\ep

\br{rfunctor}
{\em In view of Propositions \ref{p2.14} and \ref{pli4}, we obtain
a canonical functor $F_{\a}$ from the category of
weak $V$-modules to itself in the obvious way. Since
\begin{eqnarray}
\Delta(\a,z)\Delta(-\a,z)=\Delta(-\a,z)\Delta(\a,z)=1,
\end{eqnarray}
we easily see that $F_{-\a}$ is the inverse functor of $F_{\a}$. 
Therefore, $F_{\a}$ is an isomorphism.}
\er

The following result [Li4, Proposition 2.4] is a generalization
of Proposition \ref{p2.14}:

\bp{pdeformintertwiningoperator}
Let $\alpha, \Delta(\alpha,z)$ be given as in Proposition \ref{p2.14}. 
Let $I$ be an intertwining operator of type
${W_{3}\choose W_{1}W_{2}}$. Define
\begin{eqnarray}
I_{\a}(w_{1},z)w_{2}=I(\Delta(\alpha,z)w_{1},z)w_{2}
\end{eqnarray}
for $w_{1}\in W_{1},\; w_{2}\in W_{2}$.
Then $I_{\a}$ is an intertwining operator of type 
${W_{3}^{(\a)}\choose W_{1} W_{2}^{(\a)}}$.
\ep

The following results [Li4, Corollary 2.12, Theorem 2.15, Proposition 3.2]
give a construction of simple currents:

\bt{tli3} Let $V$ be a simple vertex operator algebra and 
let $\alpha\in V_{(1)}$ be such that (\ref{h2.11}) holds 
and such that $\alpha(0)$ has only integral eigenvalues on $V$.
If for each irreducible (ordinary) $V$-module $W$,
the weak module $W^{(\a)}$ is
an (ordinary) $V$-module, then $V^{(\a)}$ is a simple current.
Furthermore, for any irreducible $V$-module $W$,
\begin{eqnarray}
[W]\cdot [V^{(\a)}]=[W^{(\a)}].
\end{eqnarray}
\et

The following lemma gives more information about $V^{(\a)}$.

\bl{lweightinformation}
Let $V, \a$ be as in Theorem \ref{tli3}. Then $L(0)$ acts
semisimply on $V^{(\a)}$ with eigenvalues in ${1\over 2}\gamma+\Z$,
where $\a(1)\a=\gamma {\bf 1}$.
\el

\pf From [Li5, (3.18)], we have 
\begin{eqnarray}\label{enewvirasoro}
\Delta(\a,z)\omega=\omega+\a z^{-1}
+{1\over 2}\a(1)\a z^{-2}=\omega+\a z^{-1}
+{1\over 2}\gamma {\bf 1}z^{-2},
\end{eqnarray}
where $\omega$ is the Virasoro element. Then
\begin{eqnarray}
Y_{\a}(\omega,z)=Y(\Delta(\a,z)\omega,z)=Y(\omega,z)+z^{-1}Y(\a,z)
+{1\over 2}\gamma z^{-2}.
\end{eqnarray}
In terms of components we have
\begin{eqnarray}
L_{\a}(m)=L(m)+\a(m)+{1\over 2}\gamma \delta_{m,0}
\end{eqnarray}
for $m\in \Z$, where 
$Y_{\a}(\omega,z)=\sum_{m\in \Z}L_{\a}(m)z^{-m-2}$.
Since $L(0)$ and $\a(0)$ act semisimply on $V$ with integral
eigenvalues, $L_{\a}(0)$ acts semisimply on $V$ with 
eigenvalues in ${1\over 2}\gamma+\Z$. That is, $L(0)$ acts
semisimply on $V^{(\a)}$ with eigenvalues in ${1\over 2}\gamma+\Z$.
$\;\;\;\;\Box$

Since any irreducible weak module is an ordinary module
for a regular vertex operator algebra, from Theorem \ref{tli3}
we immediately have:

\bc{csimple}
Let $V$ be a regular vertex operator algebra and let $\alpha\in V_{(1)}$
be given as in Theorem \ref{tli3}. Then $V^{(\a)}$ is a simple current.
In particular, this is true when $V$ is a tensor product from the following algebras:
$$L_{\fg}(\ell,0),\;\;\; V_{L},\;\;\; L_{\fg}(\ell,0)\otimes V_{L},$$
where $\ell$ is a positive integer and $L$ is a positive-definite 
even lattice.$\;\;\;\;\Box$
\ec

The following result was obtained in [Li4]:

\begin{prop}\label{p2.15} 
Let $L$ be a positive definite even lattice.
(1) For $\beta\in L^{o}$, as a $V_{L}$-module,
$V_{L}^{(\b)}$ is isomorphic to $V_{L+\beta}$.
(2) Every irreducible $V_{L}$-module is a simple current.
(3) The Verlinde algebra is canonically isomorphic to the group algebra 
$\C[L^{o}/L]$.
\end{prop}

Previously,
intertwining operators for $V_{L}$ were explicitly constructed, 
fusion rules were calculated and (3) was proved in [DL]. 
(Of course, (3) implies (2).) 

Though vertex operator algebra $M_{{\bf h}}(1,0)$ is not regular,
the same proof of Proposition \ref{p2.15} gives the following result:

\begin{prop}\label{pfreefield} (1) For $h\in {\bf h}$,
as an $M_{{\bf h}}(1,0)$-module,
$M_{{\bf h}}(1,0)^{(h)}$ is isomorphic to $M_{\bf h}(1,h)$.
(2) Every irreducible $M_{{\bf h}}(1,0)$-module is a simple current.
(3) The Verlinde algebra is canonically isomorphic to the group algebra
$\C[{\bf h}]$. $\;\;\;\;\Box$
\end{prop}

\br{rtensorproduct}
{\em Suppose that $V=V^{1}\otimes V^{2}$ is a tensor product vertex 
operator algebra. Then $V^{1}_{(1)}$ and $V^{2}_{(1)}$ are
canonical subspaces of $V_{(1)}$. Suppose $\a=\a'+\a''$
where $\a^{1}\in V^{1}_{(1)},\; a^{2}\in V^{2}_{(1)}$. 
Then $\a$ satisfies (\ref{h2.11}) if and only if both $\a^{1}$ and $\a^{2}$
satisfy (\ref{h2.11}). Furthermore, $\a(0)$ acting on $V$ has only integral 
eigenvalues if and only if $\a^{i}(0)$ acting on $V^{i}$ has 
only integral eigenvalues for $i=1,2$.
Let $W=W_{1}\otimes W_{2}$, where $W_{1}, W_{2}$ are $V_{1}$ 
and $V_{2}$-modules, respectively.
Since $[\a^{1}(m),\a^{2}(n)]=0$ for $m,n\in \Z$, we have
\begin{eqnarray}
\Delta(\a,z)=\Delta(\a^{1},z)\Delta(\a^{2},z).
\end{eqnarray}
Then we have
\begin{eqnarray}
W^{(\a)}=W_{1}^{(\a^{1})}\otimes W_{2}^{(\a^{2})}.
\end{eqnarray}}
\er

Let $\fg, {\bf h}, \<\cdot,\cdot\>$ be as in Section 2.1.
Let $\{e_{i}, f_{i}\;|\;i=1,\dots, n\}$ be the Chevalley generators
with simple roots $\a_{1},\dots,\a_{n}$ and simple coroots 
$\alpha^{\vee}_{1},\dots,\alpha^{\vee}_{n}$. 
Let $\lambda_{i}$ ($i=1,\dots,n$) be the fundamental
weights for $\fg$.
Let $Q=\oplus_{i=1}^{n}\Z \a_{i}$ be the root lattice and
let $P=\oplus_{i=1}^{n}\Z \lambda_{i}$ be the weight lattice.

Let $h^{(1)},\dots, h^{(n)}\in {\bf h}$ be the fundamental 
co-weights, i.e.,
\begin{eqnarray}
\alpha_{i}(h^{(j)})=\delta_{i,j}\;\;\;\mbox{ for }i, j=1,\dots,n.
\end{eqnarray}
Set
\begin{eqnarray}
Q^{\vee}&=&\Z \a_{1}^{\vee}+\cdots +\Z \a_{n}^{\vee},\\
P^{\vee}&=&\Z h^{(1)}+\cdots +\Z h^{(n)},
\end{eqnarray}
the co-root lattice and the co-weight lattice.

Let
\begin{eqnarray}
\theta=\sum_{i=1}^{n}a_{i}\alpha_{i}
\end{eqnarray}
be the highest long root. Then
\begin{eqnarray}
\theta(h^{(i)})=a_{i}\;\;\;\mbox{ for }i=1,\dots,n.
\end{eqnarray}
We shall need to know which $a_{i}$ equal $1$. 
The following is a list for such $a_{i}$ (cf. [K]):
\begin{eqnarray}
 A_{n}:& &a_{1},\dots, a_{n}\nonumber\\
B_{n}: & &a_{1}\nonumber\\
C_{n}:& & a_{n}\nonumber\\
D_{n}: & & a_{1},a_{n-1},a_{n}\nonumber\\
E_{6}:& &a_{1},a_{5}\nonumber\\
E_{7}: & &a_{6}.
\end{eqnarray}

\br{rnumbering} {\em Simple roots for type $E$ ($E_{6}$, $E_{7}$, $E_{8}$)
were numbered differently in [H] and [K]. Here, we use
the numbering system of [K].}
\er


Let $\Lambda_{0}, \dots, \Lambda_{n}$ be the fundamental weights 
of $\hat{\fg}$ [K]. Then
each $\lambda_{i}$ for $1\le i\le n$ is naturally extended to 
$\Lambda_{i}$. From the Dynkin diagram ([K], TABLE Aff 1) we find that
$a_{i}=1$ if and only if the vertices $0$ and $i$ are in the same
orbit under the automorphism group of the affine Dynkin diagram.
We point out that if $a_{i}=1$, then $\Lambda_{i}$ is 
of level one.

The following proposition was proved in [Li4] 
(Proposition 3.5, Remark 3.8):

\begin{prop}\label{p5.1} 
Let $\ell$ be a complex number with $\ell\ne -h^{\vee}$,
where $h^{\vee}$ is the dual Coxeter number of $\fg$. If
the coefficient $a_{i}$ of $\a_{i}$ in $\theta$ is $1$, then 
as an $L(\ell,0)$-module
\begin{eqnarray}
L(\ell,0)^{(h^{(i)})}\simeq L(\ell,\ell\lambda_{i}).
\end{eqnarray}
Furthermore, if $\ell$ is a positive integer,
$L(\ell,\ell\lambda_{i})$ is a simple current for $L(\ell,0)$.
\end{prop}

\begin{rem}\label{r5.2}
{\em It was known ([FG], [F]) that $L(\ell,\ell\lambda_{i})$ with $a_{i}=1$
are all the simple currents except the level $2$ simple current 
$L(\Lambda_{7})$ for $\fg$ of type $E_{8}$.}
\end{rem}

\begin{rem}\label{r5.3} {\em The element
$\Delta(h^{(i)},z)$ gives rise to an automorphism $\psi_{i}$ of
affine Lie algebra $\hat{\fg}$ via
\begin{eqnarray}
\psi_{i}(Y(a,z))=Y_{h^{(i)}}(a,z)=Y(\Delta(h^{(i)},z)a,z)
\end{eqnarray}
for $a\in \fg$, where $Y(a,z)=\sum_{n\in \Z}a(n)z^{-n-1}$ is 
the generating series of $a$.
That is,
\begin{eqnarray}
& &\psi_{i}(\alpha^{\vee}_{i}(n))=\alpha^{\vee}_{i}(n)+\delta_{n,0}\ell,
\;\;\;\;\;
\psi_{i}(e_{i}(n))=e_{i}(n+1),\;\;\;\;\; 
\psi_{i}(f_{i}(n))=f_{i}(n-1);\\
& &\psi_{i}(\alpha^{\vee}_{j}(n))=\alpha^{\vee}_{j}(n),
\;\;\;\psi_{i}(e_{j}(n))=e_{j}(n),\;\;\;\;\; 
\psi_{i}(f_{j}(n))=f_{j}(n)
\;\;\;\mbox{for }j\ne i, n\in {\Z},\hspace{1cm}
\end{eqnarray}
and 
\begin{eqnarray}
\psi_{i} (f_{\theta}(n))=f_{\theta}(n-1)\;\;\;\mbox{for }n\in {\Z}.
\end{eqnarray}
The vector ${\bf 1}$ in 
$L(\ell,0)^{(h^{(i)})}$ is a
highest weight vector of weight $\ell\Lambda_{i}$. 
More general automorphisms of this type was recently studied in [FS].
Note that the composition of the representation on $L(\ell,0)$ with
the corresponding Dynkin diagram automorphism of $\hat{\fg}$
also gives $L(\ell,\ell\lambda_{i})$.
But $\psi_{i}$ are not Dynkin diagram automorphisms. 
Dynkin diagram automorphisms played important roles
in [FG], [F], [SY] and [FM].}
\end{rem}

\br{rlowestweight}
{\em For $\a\in {\bf h}$, from (\ref{enewvirasoro}) we have
\begin{eqnarray}
L_{\a}(m)=L(m)+\a(m)+{1\over 2}\ell \<\a,\a\>\delta_{m,0}
\end{eqnarray}
for $m\in \Z$, recall that
$Y_{\a}(\omega,z)=\sum_{m\in \Z}L_{\a}(m)z^{-m-2}$ with $\omega$ being 
the Virasoro element.
Then the $L_{h^{(i)}}(0)$-weight of ${\bf 1}$ is 
${1\over 2}\ell \<h^{(i)},h^{(i)}\>$.
Since under the new action $Y(\Delta(h^{(i)},z)\cdot,z)$ on $L(\ell,0)$,
${\bf 1}$ is still a highest weight vector for $\hat{\fg}$,
the lowest $L(0)$-weight of $L(\ell,\ell \lambda_{i})$ is
${1\over 2}\ell \<h^{(i)},h^{(i)}\>$. Of course, there is a formula
for the lowest weight of any irreducible module (cf. [DL]).}
\er

{}From now on we assume that $k$ is a positive integer. Let $h\in {\bf h}$.
Note that for any root vector $e_{\a}$ of $\fg$ and for $m\in \Z$,
\begin{eqnarray}\label{ehbracket}
[h(0),e_{\a}(m)]=\a(h)e_{\a}(m).
\end{eqnarray}
Since $h(0){\bf 1}=0$ and $L(k,0)$ is generated by
$\hat{\fg}$ from ${\bf 1}$, using (\ref{ehbracket}) we see that
$h(0)$ has only integral eigenvalues on $L(k,0)$ 
if and only if $h\in P^{\vee}$.

Define a map $\pi$ from $P^{\vee}$ to the Verlinde algebra 
${\cal{V}}(L(k,0))$ of $L(k,0)$ by
\begin{eqnarray}
\pi (\a)=[L(k,0)^{(\a)}]\;\;\;\mbox{ for }\a\in P^{\vee}.
\end{eqnarray}
The map $\pi$ naturally extends to a linear map from 
the group algebra $\C[P^{\vee}]$ to ${\cal{V}}(L(k,0))$.
We abuse the notation $\pi$ for this extension also.

\bp{phomomorphismpi}
The linear map $\pi$ is an algebra homomorphism.
\ep

\pf For $\a,\b\in P^{\vee}$, since $[\a (r),\b (s)]=0$ for $r,s\ge 0$, we have
\begin{eqnarray}
\Delta(\a+\b,z)=\Delta(\a,z)\Delta(\b,z),
\end{eqnarray}
recall (\ref{edeltadefinition}). Then we get
\begin{eqnarray}
L(k,0)^{(\a+\b)}\simeq (L(k,0)^{(\a)})^{(\b)}
\end{eqnarray}
as $L(k,0)$-modules. In view of Theorem \ref{tli3} we have
\begin{eqnarray}
[L(k,0)^{(\a)}]\cdot [L(k,0)^{(\b)}]=[(L(k,0)^{(\a)})^{(\b)}]
=[L(k,0)^{(\a+\b)}].
\end{eqnarray}
Thus $\pi$ is an algebra homomorphism. $\;\;\;\;\Box$

It follows immediately that $\pi(P^{\vee})$ is an abelian group.
The following result gives important kernel elements of $\pi$ 
as a group homomorphism on $P^{\vee}$.

\bp{pequivalence}
We have
\begin{eqnarray}\label{eequivalenceaffine}
[L(k,0)^{(\a)}]=[L(k,0)]\;\;\;\mbox{ for }\a\in Q^{\vee}.
\end{eqnarray}
Furthermore, for any irreducible $L(k,0)$-module $W$, we have
\begin{eqnarray}
[W^{(\a)}]=[W]\;\;\;\mbox{ for }\a\in Q^{\vee}.
\end{eqnarray}
\ep

\pf It suffices to prove (\ref{eequivalenceaffine}) for
$\a=\a_{i}^{\vee}$.
For $1\le i\le n$, set $r=\frac{2}{\<\a_{i},\a_{i}\>}$,
a positive integer. From [H] or [K] we have
\begin{eqnarray}
\<\a_{i}^{\vee},\a_{i}^{\vee}\>=\frac{4}{\<\a_{i},\a_{i}\>}=2r.
\end{eqnarray}
Since
\begin{eqnarray}
[e_{i}(1),f_{i}(-1)]=\a_{i}^{\vee}(0)+\<e_{i},f_{i}\>c
=\a_{i}^{\vee}(0)+{1\over 2}\<\a_{i}^{\vee},\a_{i}^{\vee}\>c
=\a_{i}^{\vee}(0)+rc,
\end{eqnarray}
${\cal{L}}_{i}:=\C e_{i}(1)+\C f_{i}(-1)+\C (\a_{i}^{\vee}+rc)$ 
is a subalgebra of 
$\hat{\fg}$ isomorphic to $sl(2)$. Furthermore, $L(k,0)$, 
being an integrable $\hat{\fg}$-module, is an integrable 
${\cal{L}}_{i}$-module.
Clearly,  ${\bf 1}$ is a highest weight vector of weight $rk$
for ${\cal{L}}_{i}$. Then $f_{i}(-1)^{rk}{\bf 1}\ne 0$.
{}From (\ref{enewvirasoro}) we have
\begin{eqnarray}
L_{\a_{i}^{\vee}}(0)
=L(0)+\a_{i}^{\vee} (0)+{1\over 2}k\<\a_{i}^{\vee},\a_{i}^{\vee}\>
=L(0)+\a_{i}^{\vee} (0)+kr,
\end{eqnarray}
where $Y_{\a_{i}^{\vee}}(\omega,z)=\sum_{n\in \Z}L_{\a_{i}^{\vee}}(n)z^{-n-2}$.
Then
\begin{eqnarray}
L_{\a_{i}^{\vee}}(0)f_{i}(-1)^{rk}{\bf 1}=(L(0)+
\a_{i}^{\vee}(0)+kr)f_{i}(-1)^{rk}{\bf 1}=(kr-2kr+kr)f_{i}(-1)^{rk}{\bf 1}=0.
\end{eqnarray}
That is, $f_{i}(-1)^{rk}{\bf 1}$ is a non-zero element of
$L(k,0)^{(\a_{i}^{\vee})}$ of weight zero.

On the other hand, with $L(k,0)$ being regular, 
every irreducible module, in particular,
 $L(k,0)^{(\a_{i}^{\vee})}$, is an integrable $\hat{\fg}$-module 
of level $k$, which is unitary. Thus
the lowest weights of $L(k,0)^{(\a_{i}^{\vee})}$ are nonnegative. 
Consequently,
the lowest weight of $L(k,0)^{(\a_{i}^{\vee})}$ is $0$.
Because $L(k,0)$ is the only irreducible $L(k,0)$-module with $0$ 
being the lowest weight,
we must have $L(k,0)^{(\a_{i}^{\vee})}\simeq L(k,0)$ as $L(k,0)$-modules.

For an irreducible $L(k,0)$-module $W$, using the first part and 
Theorem \ref{tli3} we get
\begin{eqnarray}
[W]=[W]\cdot [V]=[W]\cdot [V^{(\a)}]=[W^{(\a)}]\;\;\;\mbox{ for }
\a\in Q^{\vee}.
\end{eqnarray}
This completes the proof.$\;\;\;\;\Box$

The following result generalizes Proposition \ref{p5.1} with $L=Q$:

\bt{tsimplefusionalgebra}
The algebra homomorphism $\pi$ gives rise to an algebra isomorphism from
the group algebra $\C [P^{\vee}/Q^{\vee}]$ onto 
$\pi (\C [P^{\vee}])$. Furthermore, 
\begin{eqnarray}\label{enewfact}
\pi (P^{\vee})=\{ [L(k,k\lambda_{i})]\;|\;a_{i}=1\}.
\end{eqnarray}
\et

\pf In view of Proposition \ref{pequivalence}, $\pi$ gives rise to 
an algebra homomorphism $\bar{\pi}$ from 
$\C[P^{\vee}/Q^{\vee}]$ onto $\pi(\C [P^{\vee}])$.
{}From [H] (Section 13.1), we have
\begin{eqnarray}
|P/Q|=n+1,\; 2,\; 2,\; 4,\; 3,\; 2,\; 1,\;1,\; 1
\end{eqnarray}
for $\fg$ of type $A_{n+1}, B_{n}, C_{n}, D_{n}, E_{6}, E_{7}, E_{8}, F_{4}, 
G_{2}$, respectively. 
Note that $B_{n}$ and $C_{n}$ are dual to each other and 
the others are self-dual.
Then $|P^{\vee}/Q^{\vee}|=|P/Q|$ for all types.
On the other hand, from Proposition \ref{p5.1},
all $[L(k,k\lambda_{i})]$ for $i$ with $a_{i}=1$, which are
distinct basis elements ${\cal{V}}(L(k,0))$, are images of $\pi$.
Then it follows immediately that (\ref{enewfact}) holds and 
$\bar{\pi}$ is an algebra isomorphism.
$\;\;\;\;\Box$

The group structure of $P/Q$ for the nontrivial 
cases was given in [H] (Exercise 4 on page 71). Then we have:

For $A_{n+1}$, $P^{\vee}/Q^{\vee}\simeq \Z/(n+1)\Z$ with
$h^{(1)}+Q^{\vee}$ as a generator such that
\begin{eqnarray}
mh^{(1)}+Q^{\vee}=h^{(\bar{m})}+Q^{\vee},
\end{eqnarray}
where $\bar{m}$ 
is the least nonnegative residue of $m$ modulo $n+1$.

For $D_{n}$ with $n$ being odd, $P^{\vee}/Q^{\vee}\simeq \Z/4\Z$ with
$h^{(n)}+Q^{\vee}$ as a generator such that
\begin{eqnarray}
2h^{(n)}+Q^{\vee}=h^{(1)}+Q^{\vee},\;\;\; 3h^{(n)}+Q^{\vee}=h^{(n-1)}+Q^{\vee}.
\end{eqnarray}

For $D_{n}$ with $n$ being even, 
$P^{\vee}/Q^{\vee}\simeq \Z/2\Z \times \Z/2\Z$
with $h^{(n-1)}+Q^{\vee}$ and $h^{(n)}+Q^{\vee}$
as generators such that
\begin{eqnarray}
h^{(1)}+Q^{\vee}=h^{(n-1)}+h^{(n)}+Q^{\vee}.
\end{eqnarray}
Then with Theorem \ref{tsimplefusionalgebra} we immediately 
have the following results which have been known to physicists:

\bc{csimplecurrentfusion}
The following fusion algebra relations hold:\\
For $A_{n+1}$, for $0\le i,j\le n$,
\begin{eqnarray}
[L(k,k\lambda_{i})]\cdot [L(k,k\lambda_{j})]=[L(k,k\lambda_{\overline{i+j}})].
\end{eqnarray}
For $B_{n}$,
\begin{eqnarray}
[L(k,k\lambda_{1})]\cdot [L(k,k\lambda_{1})]=[L(k,0)].
\end{eqnarray}
For $C_{n}$,
\begin{eqnarray}
[L(k,k\lambda_{n})]\cdot [L(k,k\lambda_{n})]=[L(k,0)].
\end{eqnarray}
For $D_{n}$ with odd $n$,
\begin{eqnarray}
[L(k,k\lambda_{n})]^{2}=[L(k,k\lambda_{1})],\;\;\;
[L(k,k\lambda_{n})]^{3}=[L(k,k\lambda_{n-1})],\;\;
[L(k,k\lambda_{n})]^{4}=[L(k,0)].
\end{eqnarray}
For $D_{n}$ with even $n$,
\begin{eqnarray}
& &[L(k,k\lambda_{1})]^{2}=[L(k,k\lambda_{n-1})]^{2}=[L(k,k\lambda_{n})]^{2}
=[L(k,0)],\\
& &[L(k,k\lambda_{n-1})]\cdot [L(k,k\lambda_{n})] =[L(k,k\lambda_{1})].
\end{eqnarray}
For $E_{6}$,
\begin{eqnarray}
[L(k,k\lambda_{1})]^{2}=[L(k,k\lambda_{5})],\;\;\;
[L(k,k\lambda_{1})]^{3}=[L(k,0)].
\end{eqnarray}
For $E_{7}$,
\begin{eqnarray}
[L(k,k\lambda_{6})]^{2}=[L(k,0)].\;\;\;\;\Box
\end{eqnarray}
\ec

We shall need the number $\<h^{(i)},h^{(i)}\>$ and the explicit expression 
of $h^{(i)}$ in terms of $\a_{j}^{\vee}$.
The expression of each $\lambda_{i}$ in terms of simple roots $\a_{j}$
was known ([H], page 69, Table 1). Suppose that for $1\le i\le n$,
\begin{eqnarray}
& &\lambda_{i}=a_{i1}\a_{1}+\cdots + a_{in}\a_{n},\\
& &h^{(i)}=b_{i1}\a_{1}^{\vee}+\cdots +b_{in}\a_{n}^{\vee}.
\end{eqnarray}
Then
\begin{eqnarray}
a_{ij}=\lambda_{i}(h^{(j)})=b_{ji}.
\end{eqnarray}
Furthermore, from [H] we get
\begin{eqnarray}
\<h^{(i)}, \a_{j}^{\vee}\>=\<h^{(i)}, t_{\a_{j}}\>\frac{2}{\<\a_{j},\a_{j}\>}
=\a_{j}(h^{(i)})\frac{2}{\<\a_{j},\a_{j}\>}
=\delta_{i,j}\frac{2}{\<\a_{j},\a_{j}\>}.
\end{eqnarray}
Then
\begin{eqnarray}
\<h^{(i)},h^{(i)}\>=b_{ii}\frac{2}{\<\a_{i},\a_{i}\>}
=a_{ii}\frac{2}{\<\a_{i},\a_{i}\>}.
\end{eqnarray}
With a glance of Table Aff in [K] we see that if $a_{i}=1$, 
$\a_{i}$ is a long root, hence
$\<\a_{i},\a_{i}\>=2$.
Therefore,  we have obtained:

\bl{lsimplefact}
If $\lambda_{i}=a_{i1}\a_{1}+\cdots + a_{in}\a_{n}$ for $1\le i\le n$,
then 
\begin{eqnarray}
& &h^{(i)}=a_{1i}\a_{i}^{\vee}+\cdots+a_{ni}\a_{i}^{\vee},\\
& &\<h^{(i)},h^{(i)}\>=a_{ii}\frac{2}{\<\a_{j},\a_{j}\>}.
\end{eqnarray}
In particular, if $a_{i}=1$, we have $\<h^{(i)},h^{(i)}\>=a_{ii}$.
$\;\;\;\;\Box$
\el

\newpage

\section{Extension of vertex operator algebras by simple currents}
We shall first recall some of the results from [DLM2] on $V[L]$, an 
extension of a certain vertex operator algebra $V$ by simple currents
$V^{(\a)}$ parametrized by $\a\in L$, where $L$ is a lattice
carrying an intrinsic structure of $V$.
We then extend and refine those results. 
In Section 3.3, we concentrate a special class of $V[L]$.
We classify all irreducible $V[L]$-modules,
prove a complete reducibility theorem and we give a formula of fusion 
rules in the category of $V[L]$-modules in terms of
fusion rules in the category of $V$-modules.

\subsection{Extension of algebras}
We shall first establish some basic assumptions which 
will remain in force throughout this section.

Let $V$ be a {\em simple vertex operator algebra}.
As in Section 2, one may think of $V$ as
a tensor product from the following vertex operator algebras:
$$L_{\fg}(\ell,0),\;\;\; M_{\bf h}(1,0),\;\;\; V_{L}.$$

Let $H$ be a (finite-dimensional) subspace of $V_{(1)}$ 
satisfying the following conditions:
\begin{eqnarray}\label{e3.1}
L(n)h=\delta_{n,0}h,\;\;\;h(n)h'=B(h,h') \delta_{n,1}{\bf 1}
\;\;\;\mbox{for }n\in {\Z}_{+},\; h,h'\in H,\label{3.7}
\end{eqnarray}
where $Y(h,z)=\sum h(n)z^{-n-1}$ for $h\in H$, and
$B(\cdot,\cdot)$ is {\em assumed} to be a 
nondegenerate symmetric bilinear form on $H$. 
We then identify $H$ with its dual $H^{*}$.
We also {\em assume} that for any $h\in H$, 
$h(0)$ acts semisimply on $V$. Then
\begin{eqnarray}
V=\oplus_{\alpha\in H}V^{(0,\alpha)},\;\;\;\mbox{ where }
V^{(0,\alpha)}=\{v\in V\;| \;h(0)v=B(\alpha,h) v
\;\;\mbox{for }h\in H\}.
\end{eqnarray}
Set 
\begin{eqnarray}
P=\{\a \in H\;|\;V^{(0,\a)}\ne 0\}.
\end{eqnarray}
As $V$ is simple, $P$ is a subgroup of $H$ (cf. [LX]).
Then $P$ equipped with the bilinear form $B$ is a lattice.

Let $L$ be a subgroup of $H$ such that for each $\alpha\in L$, 
$\alpha(0)$ acting on $V$ has only {\em integral eigenvalues}. 
This amounts to $L\subset P^{o}$, where
\begin{eqnarray}
P^{o}=\{h\in H\;|\; B(h,\a)\in \Z\;\;\mbox{ for }\a\in P\}
\end{eqnarray}
is the dual lattice of $P$. Note that the rank of $P$ may be less than 
$\dim H$.

Let $W$ be a $V$-module. By Proposition \ref{p2.14},
for $\alpha\in L$, we have a (weak) $V$-module 
\begin{eqnarray}
(W^{(\alpha)}, Y_{\alpha}(\cdot,z)):=(W,Y(\Delta(\alpha,z)\cdot,z)).
\end{eqnarray}
For convenience, we reformulate the construction of the $V$-module
$W^{(a)}$ as follows: Set
\begin{eqnarray}
W^{(\a)}=\C e^{\a}\otimes W\simeq W\;\;(\mbox{linearly}),
\end{eqnarray}
where $e^{\a}$ is a symbol for now and $\C e^{\a}$ is a one-dimensional 
vector space with $e^{\a}$ as a pre-chosen basis element. Then define
\begin{eqnarray}
Y_{\a}(v,z)(e^{\a}\otimes w)=e^{\a}\otimes Y(\Delta(\a,z)v,z)w
\;\;\;\mbox{ for }v\in V,\; w\in W.
\end{eqnarray}

Set 
\begin{eqnarray}
W[L]=\oplus _{\alpha\in L}W^{(\alpha)}=\C[L]\otimes W,
\end{eqnarray}
equipped with the direct sum $V$-module structure. Now, 
the symbol $e^{\a}$ in the definition of $W^{(\a)}$ is considered 
as an element of the group algebra $\C[L]$.

For $\a\in L$, we define
a linear endomorphism $\psi_{\alpha}$
\footnote{The map $\psi_{\a}$ here is
the inverse of the map $\psi_{\a}$ in [DLM2]} 
of $W[L]$ by
\begin{eqnarray}
\psi_{\a}(e^{\b}\otimes w)=e^{\a+\b}\otimes w\;\;\;\mbox{ for }
\b\in L,\; w\in W.
\end{eqnarray}
Then
\begin{eqnarray}
\psi_{0}=1,\;\;\;\psi_{\alpha+\beta}=\psi_{\alpha}\psi_{\beta}
\;\;\;\mbox{ for }\alpha, \beta\in L.
\end{eqnarray}
That is, $\psi$ gives rise to a representation of $L$ on $W[L].$ 

For $\a\in L$, we set ([LW], [FLM])
\begin{eqnarray}\label{eepmaz}
E^{\pm}(\a,z)
=\exp \left(\sum_{n=1}^{\infty}\frac{h(\pm n)}{\pm n}z^{\mp n}\right).
\end{eqnarray}
Then
\begin{eqnarray}
\Delta(\a,z)=z^{\a(0)}E^{+}(-\a,-z).
\end{eqnarray}

Next, we extend the domain of $Y_{\a}$ from $V$ to $V[L]$.

\bd{d3.5} {\em For $u\in V^{(\alpha)},\; v\in V^{(\beta)}$ with 
$\alpha,\beta\in L$, we define 
\be{edefvoa}
Y_{\a}(u,z)v=
\psi_{\alpha+\beta}E^{-}(-\alpha,z)Y(\psi_{-\alpha}\Delta(\beta,z)u,z)
\Delta(\alpha,-z)\psi_{-\beta}(v)\in V^{(\alpha+\beta)}\{z\}.
\ee
We then define a linear map $\tilde{Y}(\cdot,z)$ from 
$V[L]$ to $(\End V[L])\{z\}$ via $\tilde{Y}(u,z)=Y_{\a}(u,z)$ for
$u\in V^{(\a)}$.} 
\ed

Note that the $E^{-}(\a,z)$ defined in [DLM2] is
the $E^{-}(-\a,z)$ defined in (\ref{eepmaz}) ([LW], [FLM]).
Then the definition of $Y_{\a}$ is is exactly the same as the 
one defined in [DLM2].


For $\alpha\in P,\; h\in H$, set
\begin{eqnarray}
V^{(\alpha,h)}=\psi_{\a}(V^{(0,h)}).
\end{eqnarray}
Then 
\begin{eqnarray}
V^{(\alpha)}=\oplus_{h\in P}V^{(\alpha,h)}.
\end{eqnarray}

\bd{dceta}
{\em We define $\C$-valued functions $\eta$ and $C$ on 
$(L\times H)\times (L\times H)$ by
\be{3.19}
\eta((\alpha_{1},h_{1}),(\alpha_{2},h_{2}))=
-B(\alpha_{1},\a_{2})
-B(\alpha_{1},h_{2})-B(\alpha_{2},h_{1})\in \C,
\ee
\be{3.20}
C((\alpha_{1},h_{1}),(\alpha_{2},h_{2}))=
e^{(B(\alpha_{1},h_{2})-B(\alpha_{2},h_{1}))\pi i}\in {\C}^{\times}
\ee
for $(\alpha_{i},h_{i})\in L\times H,\; i=1,2$.}
\ed

Then we have ([DLM2], Theorem 3.5):

\bt{t3.7} Let $u\in V^{(\alpha,h_{1})},
\; v\in V^{(\beta,h_{2})}, \; w\in V^{(\gamma,h_{3})}$ with
$\alpha,\beta,\gamma\in L,\; h_{1},h_{2},h_{3}\in P$. Then
\begin{eqnarray}
& &z_{0}^{-1}\delta\left(\frac{z_{1}-z_{2}}{z_{0}}\right)
\left(\frac{z_{1}-z_{2}}{z_{0}}\right)^{\eta((\alpha,h_{1}),(\beta,h_{2}))}
\tilde{Y}(u,z_{1})\tilde{Y}(v,z_{2})w\nonumber\\
&-&C((\alpha,h_{1}),(\beta,h_{2}))z_{0}^{-1}\delta\left(\frac{z_{2}-z_{1}}
{-z_{0}}\right)
\left(\frac{z_{2}-z_{1}}{z_{0}}\right)^{\eta((\alpha,h_{1}),(\beta,h_{2}))}
\tilde{Y}(v,z_{2})\tilde{Y}(u,z_{1})w\nonumber\\
&=&z_{2}^{-1}\delta\left(\frac{z_{1}-z_{0}}{z_{2}}\right)
\left(\frac{z_{2}+z_{0}}{z_{1}}\right)^{\eta((\alpha,h_{1}),(\gamma,h_{3}))}
\tilde{Y}(\tilde{Y}(u,z_{0})v,z_{2})w.\label{3.31}
\end{eqnarray}
Furthermore, for all $v\in V[L]$,
\begin{eqnarray}
[L(-1),\tilde{Y}(v,z)]={d\over dz}\tilde{Y}(v,z).
\end{eqnarray}
\et

Now we express $\tilde{Y}$ more explicitly.

\bl{lexplicitya}
For $\a,\b\in L,\; u,v\in V$, we have
\begin{eqnarray}\label{eexplicitformula1}
\tilde{Y}(e^{\a}\otimes u,z)(e^{\b}\otimes v)=e^{\a+\b}\otimes
z^{B(\a,\b)}E^{-}(-\a,z)Y(\Delta(\b,z)u,z)E^{+}(-\a,z)(-z)^{\a(0)}v.
\end{eqnarray}
In particular,
\begin{eqnarray}
\tilde{Y}(e^{\a}\otimes {\bf 1},z)(e^{\b}\otimes v)=e^{\a+\b}\otimes
z^{B(\a,\b)}E^{-}(-\a,z)E^{+}(-\a,z)(-z)^{\a(0)}v.
\end{eqnarray}
\el

\pf {}From Lemma 3.2 of [DLM2] we have
\begin{eqnarray}
\psi_{-\a}\Delta(\b,z)=z^{B(\a,\b)}\Delta(\b,z)\psi_{-a}.
\end{eqnarray}
Note that the map $\psi_{\a}$ defined here is the map $\psi_{-\a}$
defined in [DLM2]. Then (\ref{eexplicitformula1}) follows immediately.
$\;\;\;\;\Box$

\br{rcompare1}
{\em Let $G=L\times P$ be the product group.   
Suppose that there exists a positive integer $T$ such that
$\eta$ restricted on $G$ is ${1\over T}{\Z}/2{\Z}$-valued.
The original theorem states that 
$(V[L], {\bf 1},\omega, \tilde{Y}, T, G, \eta(\cdot,\cdot),
C(\cdot,\cdot))$ is a 
generalized vertex algebra in the sense of [DL].
This result is similar to Theorem 5.1 of [DL], which states that
if $L$ is a rational lattice, then $V_{L}$ has a canonical
generalized vertex algebra structure. In fact, by taking $V=M_{\bf h}(1,0)$
with ${\bf h}=\C\otimes_{\Z}L$, we have $P=0$ and $V[L]=V_{L}$.}
\er

In order to have vertex (super)algebras $V[L]$, we shall restrict ourselves to 
special $L$. We have already assumed that $L\subset P^{o}$, 
or what is equivalent to, $\a(0)$ acting on $V$ has only integral 
eigenvalues for every $\a\in L$.
Now we furthermore assume that $(L,B)$ is an integral lattice. Then
\begin{eqnarray}
B(\a,\a),\;\;\;\;B(\a,\beta) \in {\Z}\;\;\;
\mbox{ for }\a\in L,\;\beta\in P.
\end{eqnarray}
Recall from Lemma \ref{lweightinformation} that
for $\a\in L$, the weights of $V^{(\a)}$ are
contained in ${1\over 2}B(\a,\a)+\Z$. Then
$V[L]$ is ${1\over 2}\Z$-graded by $L(0)$-weights.
Furthermore, the function $\eta$ restricted to 
$(L\times P)\times (L\times P)$ is $\Z$-valued
and
\begin{eqnarray}
C((\a_{1},h_{1}),(\a_{2},h_{2}))=(-1)^{B(\a_{1},h_{2})-B(\a_{2},h_{1})},
\end{eqnarray}
(recall (\ref{3.19}) and (\ref{3.20})).
Then we have ([DLM2], (3.59)):

\bc{c3.13}  
Assume that $L\subset P^{0}$ and $(L,B)$ is an integral lattice.
For $a\in V^{(\alpha)},\;b\in V^{(\beta)}$ 
with $\alpha,\beta\in L$, we have
\begin{eqnarray}\label{esuperjacobi}
& &z_{0}^{-1}\delta\left(\frac{z_{1}-z_{2}}{z_{0}}\right)
\tilde{Y}(a,z_{1})\tilde{Y}(b,z_{2})
-(-1)^{B(\a,\beta)}z_{0}^{-1}\delta\left(\frac{z_{2}-z_{1}}
{-z_{0}}\right)
\tilde{Y}(b,z_{2})\tilde{Y}(a,z_{1})\nonumber\\
& &=z_{2}^{-1}\delta\left(\frac{z_{1}-z_{0}}{z_{2}}\right)
\tilde{Y}(\tilde{Y}(a,z_{0})b,z_{2}).
\end{eqnarray}
\ec

\br{rcomment1}
{\em Set
\begin{eqnarray}
L^{e}=\{\a\in L\;|\; B(\a,\a)\in 2\Z\}.
\end{eqnarray}
Clearly, $L^{e}$ is a subgroup of $L$ of index 2.
If $L^{e}\subset 2L^{o}$, i.e.,
\begin{eqnarray}\label{eextra}
B(\a,\beta)\in 2\Z\;\;\;\mbox{ for }\a\in L^{e},\;\beta\in L,
\end{eqnarray}
we easily see that $V[L]=\oplus_{\beta\in L}V^{(\beta)}$ is a 
vertex superalgebra with
\begin{eqnarray}\label{eparity}
V[L]^{0}=V[L^{e}]=\oplus_{\a\in L^{e}}V^{(\a)},\;\;\;
V[L]^{1}=\oplus_{\a\in L-L^{e}}V^{(\a)}.
\end{eqnarray}
In particular, if $L$ is of rank one, clearly (\ref{eextra}) holds, 
hence $V[L]$ is a vertex (super)algebra.
However, without assuming (\ref{eextra})
$V[L]$ equipped with the vertex operator map $\tilde{Y}$
may not be a vertex superalgebra.
Even if $L$ is even, $V[L]$ may not be a vertex algebra
unless $B$ is $2\Z$-valued. 
Corollary 3.13 of [DLM2], which
states that $V[L]$ equipped with $\tilde{Y}$ is a vertex superalgebra 
if $L$ is integral (without the condition (\ref{eextra}),
is incorrect.}
\er

As in [FLM] and [DL] for $V_{L}$,
we shall need a $2$-cocycle $\epsilon$ on $L$.
Suppose that $(L,B)$ is an integral lattice of finite rank.
Let $\{\a_{1},\dots, \a_{d}\}$ be a $\Z$-basis for $L$.
Define $\epsilon$ to be the (uniquely determined)
$\{\pm 1\}$-valued multiplicative function
on $L\times L$ such that
\begin{eqnarray}\label{e3.30}
\epsilon(\a_{i},\a_{j})
=(-1)^{B(\a_{i},\a_{j})+B(\a_{i},\a_{i})B(\a_{j},\a_{j})}
\;\;\;\mbox{ if }i\ge j,\;
\mbox{ and }1\;\;\mbox{ otherwise}
\end{eqnarray}
(cf. [DL], [FLM]).
Note that $\e(\a_{i},\a_{i})=1$ because 
$B(\a_{i},\a_{i})+B(\a_{i},\a_{i})B(\a_{i},\a_{i})$ is even.
Then 
\begin{eqnarray}
\epsilon(\a,\beta)\epsilon(\beta,\a)^{-1}=(-1)^{B(\a,\b)+B(\a,\a)B(\b,\b)}
\;\;\;\mbox{ for }\a,\beta\in L.
\end{eqnarray}
Define a vertex operator map $Y$ on $V[L]$ by
\begin{eqnarray}
Y(e^{\a}\otimes u,z)(e^{\beta}\otimes v)
=\epsilon(\a,\b)\tilde{Y}(e^{\a}\otimes u,z)(e^{\beta}\otimes v)
\end{eqnarray}
for $\a,\b\in L,\; u,v\in V$. By Lemma \ref{lexplicitya},
\begin{eqnarray}\label{edefnewy}
& &Y(e^{\a}\otimes u,z)(e^{\beta}\otimes v)\nonumber\\
&=&\e(\a,\b)e^{\a+\b}\otimes
z^{B(\a,\b)}E^{-}(\a,z)Y(\Delta(\b,z)u,z)E^{+}(-\a,z)(-z)^{\a(0)}v.
\end{eqnarray}
{}From Corollary \ref{c3.13} we immediately obtain:

\bp{pvoa}
Suppose that $L\subset P^{o}$ and that
$(L,B)$ is an integral lattice of finite rank. Let
$\epsilon$ be the $\{\pm 1\}$-valued multiplicative function
on $L\times L$ defined in (\ref{e3.30}).
Then $V[L]$ equipped with the vertex operator map $Y$ 
defined in (\ref{edefnewy}) is 
a vertex (super)algebra
with the even and odd parts being defined in (\ref{eparity}).
$\;\;\;\;\Box$
\ep

\br{rcompare2}
{\em It was proved in [DL] (Theorem 6.7 and Remarks 6.17 and 12.38)
that if $L$ is an integral lattice, $V_{L}$ 
is a vertex superalgebra.}
\er

Since we in this paper are mainly interested in vertex operator 
(super)algebras,
for the rest of Section 3 we shall assume 
that $L\subset P^{o}$ and $(L,B)$ is integral, and 
we fix the multiplicative function $\e$.


\subsection{Extensions of modules and intertwining operators}
We continue Section 3.1 to study the extension of
an irreducible $V$-module, following [DLM2]. 
The extension $W[L]$ of an irreducible $V$-module $W$
is in general a twisted $V[L]$-module with respect to an automorphism
of $V[L]$. We shall also study
the extension of an intertwining operator.

Let $W$ be an irreducible $V$-module.
By definition, homogeneous subspaces of $W$ are finite-dimensional. 
Since $[L(0), h(0)]=0$ for $h\in H$, $H$ preserves each 
homogeneous subspace of $W$ 
so that there exist $0\ne w\in W,\; h\in H=H^{*}$ such that 
$h(0)w=B(h,h')w$ for $h'\in H$. Since $H$ acts semisimply on $V$ 
(by assumption) and
$w$ generates $W$ by $V$ (from the irreducibility of $W$), $H$ also acts
semisimply on $W$.
For any $h\in H$, we define
\begin{eqnarray}
W^{(0,h)}=\{w\in W\;|\;h(0)w=B(h,h')w\;\;\mbox{for }h'\in H\}.
\end{eqnarray}
Set
\begin{eqnarray}
P(W)=\{h \in H\;|\;W^{(0,h)}\ne 0\}.
\end{eqnarray}
Since $W$ is irreducible, $P(W)$ is an irreducible $P(V)$-set. 
Then $P(W)=h+P(V)$ for any $h\in P(W)$.

\bd{ddefinitionsigma}
{\em Let $W$ be an irreducible $V$-module and $h\in P(W)$.
Define a linear endomorphism $\sigma_{W}$ of $V[L]$ by
\begin{eqnarray}
\sigma_{W}(a)=e^{-2\pi i B(\a,h)}a\;\;\;\mbox{ for }a\in V^{(\a)}
\;\;\mbox{ with }\;\a\in L.
\end{eqnarray}}
\ed

Because $L\subset P(V)^{o}$ and $P(W)=h+P(V)$, $\sigma_{W}$ 
is well defined, i.e., it does not depend on the choice of $h$.

\bl{lautomorphism}
The defined linear endomorphism $\sigma_{W}$ of $V[L]$ is an automorphism 
of the vertex (super)algebra and 
 $\sigma_{W}=1$ if and only if $\a(0)$ has only integral eigenvalues on $W$ 
for every $\a\in L$, i.e., $P(W)\subset L^{o}$.
Furthermore, if $V$ is finitely generated, $\sigma_{W}$ is of finite order
if and only if $\alpha(0)$ has rational eigenvalues on $W$ 
for every $\alpha\in L$,  or equivalently, $B(\a,h)\in \Q$ for $\a\in L$.
\el

\pf In view of (\ref{edefvoa}),
clearly, $\sigma_{W}$ is an automorphism of the vertex (super)algebra and
$\sigma_{W}=1$ if and only if $\a(0)$ has only integral eigenvalues on $W$ 
for every $\a\in L$.
Furthermore, when $V$ is finitely generated, $\sigma_{W}$ is of finite order
if and only if $\alpha(0)$ has rational eigenvalues on $W$ 
for every $\alpha\in L$,  or equivalently, $B(\a,h)\in \Q$ for $\a\in L$. 
$\;\;\;\;\Box$

Recall that $W[L]=\C[L]\otimes W=\oplus _{\alpha\in L}W^{(\alpha)}$.

\bd{dextensionmodule}
{\em Let $W$ be an irreducible $V$-module.
For $a\in V^{(\alpha)},\; w\in W^{(\beta)},\; \alpha,\beta\in L$, 
we define $Y_{W[L]}(a,z)w\in W^{(\alpha+\beta)}\{z\}$ by}
\begin{eqnarray}\label{edefmodule}
Y_{W[L]}(a,z)w
=\e(\a,\b)\psi_{\alpha+\beta}E^{-}(\alpha,z)Y(\psi_{-\alpha}
\Delta(\beta,z)a,z)\Delta(\alpha,-z)\psi_{-\beta}(w)
\end{eqnarray}
\ed
(cf. (\ref{edefvoa})).

In terms of the notion of twisted module as defined in
[Le], [D2] and [FFR] we have ([DLM2], Theorem 3.6, Corollary 3.14):

\bp{p3.14}  Let $W$ be an irreducible $V$-module
such that $\alpha(0)$ has rational eigenvalues on $W$ 
for every $\alpha\in L$. Then the following
twisted Jacobi identity holds for $u\in V^{(\a)},\; v\in V^{(\beta)},\; 
w\in W^{(\gamma,h)}$,
\begin{eqnarray}
& &z^{-1}_0\delta\left(\frac{z_1-z_2}{z_0}\right)
Y(u,z_1)Y(v,z_2)w
-(-1)^{B(\a,\a)B(\b,\b)}z^{-1}_0\delta\left(\frac{z_2-z_1}{-z_0}\right)
Y(v,z_2)Y(u,z_1)w\nonumber\\
&=&z_2^{-1}\left(\frac{z_1-z_0}{z_2}\right)^{B(\a,h)}
\delta\left(\frac{z_1-z_0}{z_2}\right)
Y(Y(u,z_0)v,z_2)w.
\end{eqnarray}
Moreover, $W[L]$ is a $\sigma_{W}$-twisted $V[L]$-module. In particular,
if $P(W)\subset L^{o}$, i.e., $\a(0)$ acting on $W$ has only integral 
eigenvalues for $\a\in L$, $W[L]$ is an (untwisted) $V[L]$-module.
\ep

Next, we prove a functorial property.

\bp{pextensionhom2}
Let $\sigma$ be an automorphism of $V[L]$ of finite order such that
\begin{eqnarray}
& &\sigma(V^{(\a)})=V^{(\a)}\;\;\;\mbox{ for }\a\in L,\\
& &\sigma(v)=v\;\;\;\mbox{ for }v\in V=V^{(0)}.
\end{eqnarray}
Let $M$ be a $\sigma$-twisted weak $V[L]$-module which is a
direct sum of irreducible $V$-modules and
let $W$ be an irreducible $V$-submodule of $M$.
Then $\sigma=\sigma_{W}$ and
there is a canonical $V[L]$-homomorphism $\phi[L]$ from 
$W[L]$ to $M$ extending the embedding $\phi$ of $W$ into $M$. 
\ep

\pf The following arguments are essentially
the ones of [DLM2], Corollary 3.15  and Lemma 4.3.

Note that for $\alpha\in L$, $Y_{W[L]}$ restricted to $V^{(\a)}\times W$
 is a nonzero intertwining operator
of type ${W^{(\a)}\choose V^{(\a)}W}$.
Since $Y_{M}(\cdot,z)\phi$ restricted to $V^{(\a)}\times W$ 
is an intertwining operator of type
${M\choose V^{(\a)}W}$ and $[V^{(\a)}]\cdot [W]=[W^{(\a)}]$, 
it follows from Schur lemma (cf. [FHL]) that there exists a 
unique $V$-homomorphism 
$\phi_{\a}$ from $W^{(\a)}$ to $M$ such that
\begin{eqnarray}\label{e3.52}
Y_{M}(a,z)w=\phi_{\a}(Y_{W[L]}(a,z)w)
\end{eqnarray}
for $a\in V^{(\a)},\; w\in W$. 
{}From the definition of a twisted module we have
\begin{eqnarray}
z^{r_{\a}}Y_{M}(a,z)w\in M((z)),\;\; 
\phi_{\a}(z^{B(\a,h)}Y_{W[L]}(a,z)w)\in M((z)).
\end{eqnarray}
Then it follows from (\ref{e3.52}) that $r_{\a}-B(\a,h)\in \Z$, 
hence $\sigma(a)=\sigma_{W}(a)$.
Thus $\sigma=\sigma_{W}$.

Define a 
$V$-homomorphism $\phi[L]$ from $W[L]$ to $M$ by
$\phi[L]=\phi_{\a}$ on $W^{(\a)}$ for $\a\in L$.
Now we show that $\phi[L]$ is a $V[L]$-homomorphism.

Let $w\in W,\; a\in V^{(\a)},\; b\in V^{(\beta)}$ with $\a,\beta\in L$.
Let $k_{0}$ be a positive integer such that
\begin{eqnarray}
& &(z_{0}+z_{2})^{k_{0}+B(\a,h)}Y_{W[L]}(a,z_{0}+z_{2})Y_{W[L]}(b,z_{2})w
\nonumber\\
&=&(z_{2}+z_{0})^{k_{0}+B(\a,h)}Y_{W[L]}(Y(a,z_{0})b,z_{2})w,\\
& &(z_{0}+z_{2})^{k_{0}+B(\a,h)}Y_{M}(a,z_{0}+z_{2})Y_{M}(b,z_{2})\phi(w)
\nonumber\\
&=&(z_{2}+z_{0})^{k_{0}+B(\a,h)}Y_{M}(Y(a,z_{0})b,z_{2})\phi(w).
\label{e3.54}
\end{eqnarray}
Then using (\ref{e3.52}) we get
\begin{eqnarray}
& &(z_{0}+z_{2})^{k_{0}+B(\a,h)}\phi_{\a+\beta}
Y_{W[L]}(a,z_{0}+z_{2})Y_{W[L]}(b,z_{2})w
\nonumber\\
&=&(z_{2}+z_{0})^{k_{0}+B(\a,h)}\phi_{\a+\beta}Y_{W[L]}(Y(a,z_{0})b,z_{2})w
\nonumber\\
&=&(z_{2}+z_{0})^{k_{0}+B(\a,h)}Y_{M}(Y(a,z_{0})b,z_{2})\phi(w)\nonumber\\
&=&(z_{0}+z_{2})^{k_{0}+B(\a,h)}Y_{M}(a,z_{0}+z_{2})Y_{M}(b,z_{2})\phi(w)
\nonumber\\
&=&(z_{0}+z_{2})^{k_{0}+B(\a,h)}Y_{M}(a,z_{0}+z_{2})\phi_{\b}Y_{M}(b,z_{2})w.
\end{eqnarray}
Multiplying both sides by $(z_{0}+z_{2})^{-k_{0}-B(\a,h)}$ we get
\begin{eqnarray}
\phi_{\a+\beta}Y_{W[L]}(a,z_{0}+z_{2})Y_{W[L]}(b,z_{2})w
=Y_{M}(a,z_{0}+z_{2})\phi_{\b}Y_{M}(b,z_{2})w,
\end{eqnarray}
that is,
\begin{eqnarray}
\phi[L]Y(a,z_{0}+z_{2})Y_{W[L]}(b,z_{2})w
=Y_{M}(a,z_{0}+z_{2})\phi[L]Y_{M}(b,z_{2})w.
\end{eqnarray}
Since $W^{(\beta)}$ is linearly spanned by $b_{n}W$ 
for $b\in V^{(\beta)},\; n\in \Z$, we have
\begin{eqnarray}
\phi[L] (Y_{W[L]}(a,z)u)=Y_{M}(a,z)u
\end{eqnarray}
for $a\in V^{(\a)},\; u\in W^{(\beta)}$.
Thus $\phi[L]$ is a $V[L]$-homomorphism. $\;\;\;\;\Box$

Our next result gives a characterization of the equivalence relation 
on (twisted) $V[L]$-modules $W[L]$ in terms of the equivalence 
of $V$-modules:

\bp{pequivalencew[l]}
Let $W_{1}$ and $W_{2}$ be irreducible $V$-modules on which 
$\a(0)$ has rational eigenvalues for each $\a\in L$.
Then $\sigma_{W_{1}}=\sigma_{W_{2}}$ and $W_{1}[L]\simeq W_{2}[L]$
if and only if $W_{2}\simeq W_{1}^{(\a)}$ for some $\a\in L$.
\ep

\pf The ``only if'' part is clear. 
Note that
$\sigma_{W^{(\a)}}=\sigma_{W}$ for any irreducible 
$V$-module $W$ and $\a\in L$ because
$P(W^{(\a)})=\a+P(W)$. 
Assume  $W_{2}\simeq W_{1}^{(\a)}$ for some $\a\in L$. Then
$\sigma_{W_{1}}=\sigma_{W_{2}}$.
Let $\phi$ be a $V$-isomorphism from $W_{2}$ to $W_{1}^{(\a)}\subset W_{1}[L]$.
It follows from Proposition \ref{pextensionhom2} that
$\phi$ extends to a $V[L]$-homomorphism $\phi[L]$ from $W_{2}[L]$ into
$W_{1}[L]$ with $\phi[L](W_{2}^{(\beta)})= W_{1}^{(\a+\beta)}$ 
for $\beta\in L$.
With each $W^{(\b)}$ being an irreducible $V$-module, 
$\phi[L]$ is an isomorphism. $\;\;\;\;\Box$

Next, we shall extend an intertwining operator $I$ 
in the category of $V$-modules to an
intertwining operator $I[L]$ in the category of $V[L]$-modules.

\bd{dextensionintertwining}
{\em Let $W_{1},W_{2}$ and $W_{3}$ be irreducible
$V$-modules and $I$ be an intertwining operator of type 
${W_{3}\choose W_{1}W_{2}}$. We define a linear map
\begin{eqnarray}
I[L]: W_{1}[L]\rightarrow (\Hom (W_{2}[L],W_{3}[L]))\{z\}
\end{eqnarray}
by
\begin{eqnarray}
I[L](a,z)w
=\e(\a,\b)\psi_{\alpha+\beta}E^{-}(\alpha,z)I(\psi_{-\alpha}\Delta(\beta,z)a,z)
\Delta(\alpha,-z)\psi_{-\beta}(w)
\end{eqnarray}
(cf. (\ref{edefvoa}) and (\ref{edefmodule}))
for $a\in W_{1}^{(\alpha)},\; w\in W_{2}^{(\beta)}$ 
with $\alpha,\beta\in L$}.
\ed

The same proof of Theorem 3.5 in [DLM2] gives:
\begin{eqnarray}
I[L](L(-1)a,z)={d\over dz}I[L](a,z)\;\;\;\mbox{ for }a\in W_{1}[L],
\end{eqnarray}
and
\begin{eqnarray}
& &z_{0}^{-1}\delta\left(\frac{z_{1}-z_{2}}{z_{0}}\right)
\left(\frac{z_{1}-z_{2}}{z_{0}}\right)
^{\eta((\alpha,h),(\beta,h_{1}))}
Y_{W}(a,z_{1})I[L](b,z_{2})u\nonumber\\
&-&C((\alpha,h),(\beta,h_{1}))
z_{0}^{-1}\delta\left(\frac{z_{2}-z_{1}}
{-z_{0}}\right)
\left(\frac{z_{2}-z_{1}}{z_{0}}\right)
^{\eta((\alpha,h),(\beta,h_{1}))}
I[L](b,z_{2})Y_{W}(a,z_{1})u\nonumber\\
& &=z_{2}^{-1}\delta\left(\frac{z_{1}-z_{0}}{z_{2}}\right)
\left(\frac{z_{2}+z_{0}}{z_{1}}\right)
^{\eta((\alpha,h),(\gamma,h_{2}))}
I[L](Y(a,z_{0})b,z_{2})u
\end{eqnarray}
for
$a\in V^{(\alpha,h)},\;b\in W_{1}^{(\beta,h_{1})},\;
u\in W_{2}^{(\gamma,h_{2})}$, 
with $\alpha,\beta,\gamma\in L,\;h\in P,\;h_{1}\in P(W_{1}),
\;h_{2}\in P(W_{2})$. 
If the extensions $W_{i}[L]$
are (untwisted) $V[L]$-modules, then $P(W_{i})\subset P^{o}$,
so that $\eta$ and $C$ have integer values. 
Then we conclude:

\bp{pextensionintertwining}
Let $W_{1},W_{2}, W_{3}$ be irreducible $V$-modules such that
for $\a\in L$, $\a(0)$ has only integral eigenvalues on $W_{i}$
for $i=1,2,3$, or what is equivalent to, the extensions $W_{i}[L]$
are (untwisted) $V[L]$-modules. Let $I$ be an intertwining 
operator of type ${W_{3}\choose W_{1}W_{2}}$ in the category of 
$V$-modules. Then $I[L]$ is an intertwining operator of type
${W_{3}[L]\choose W_{1}[L]W_{2}[L]}$ in the category of $V[L]$-modules.
$\;\;\;\;\Box$
\ep

Note that various $V$-submodules $V^{(\alpha)}$ of
$V[L]$ may be $V$-isomorphic to each other. 
It was proved in [DLM2] that $V[L]$ contains an ideal $I$ 
such that each irreducible $V$-module $V^{(\alpha)}$ is 
of multiplicity-one in the quotient algebra $V[L]$ modulo $I$.
In the following we present an abstract reformulation of this result.

Consider an (abstract) vertex operator (super)algebra 
$U=\oplus_{g\in G}V^{g}$ graded by a (finite or infinite) abelian group $G$
satisfying the following conditions:

(C1) $V^{0}$ is a vertex operator subalgebra and 
    $V^{g}$ for $g\in G$ are simple currents for $V^{0}$.

(C2) For $u\in V^{g}, v\in V^{h}$ with $g,h\in G$, $u_{n}v\in V^{g+h}$ 
    for $n\in \Z$.

(C3) For $g,h\in G$, $u_{n}v\ne 0$ for some
    $u\in V^{g},\; v\in V^{h},\;n\in \Z$.

It is easy to see that under
these conditions, $U$ is a simple $G$-graded algebra, 
i.e., there is no nontrivial $G$-graded ideal. 
{}From Conditions (2) and (3), $Y$ restricted to $V^{g}\times V^{h}$
is a nonzero intertwining operator of type 
${V^{g+h}\choose V^{g}V^{h}}$. Then from Condition (1)
we have
\begin{eqnarray}\label{efusionformula}
[V^{g}]\cdot [V^{h}]=[V^{g+h}]\;\;\;\mbox{ for }g,h\in G.
\end{eqnarray}

Set
\begin{eqnarray}
G_{0}=\{g\in G\;|\; V^{g}\simeq V^{0}\;\;\;\mbox{ as }V^{0}\mbox{-modules }\}.
\end{eqnarray}
Using (\ref{efusionformula}) by a routine argument
we easily get (cf. [DLM2], Lemma 3.7):

\bl{labstract1}
The defined subset $G_{0}$ of $G$ is a subgroup and for $g,h\in G$,
$[V^{g}]=[V^{h}]$ if and only if $g-h\in G_{0}$.$\;\;\;\;\Box$
\el


For each $g_{0}\in G_{0}$, fix  a 
$V^{0}$-isomorphism $f_{g_{0}}$ from $V^{0}$ to $V^{g_{0}}$.
We particularly define $f_{0}=1$.
Let $g_{0}\in G_{0},\; h\in G$. Then $Y(\cdot,z)\circ f_{g_{0}}$ 
is a nonzero intertwining
operator of type ${V^{g+h}\choose V^{h}V^{0}}$.
On the other hand,
for any $V^{0}$-isomorphism $\psi$ from $V^{h}$ to $V^{g+h}$,
$\psi \circ Y(\cdot,z)$ is a nonzero intertwining operator
of type ${V^{g+h}\choose V^{h}V^{0}}$.
Because $V^{h}, V^{0}, V^{g+h}$ are simple currents and
$$[V^{g+h}]=[V^{g}]\cdot [V^{h}]=[V^{0}]\cdot [V^{h}]=[V^{h}],$$
there exists a unique $V^{0}$-isomorphism $f_{g_{0},h}$ from 
$V^{h}$ to $V^{g_{0}+h}$ such that
\begin{eqnarray}
f_{g_{0},h}(Y(u,z)v)=Y(u,z)f_{g_{0}}(v)\;\;\;
\mbox{ for }v\in V^{0},\; u\in V^{h}.
\end{eqnarray}
Define a $V^{0}$-endomorphism $\bar{f}_{g_{0}}$ of $U$ via
$\bar{f}_{g_{0}}=f_{g_{0},h}$ on $V^{h}$ for $h\in G$.

Next we define $I$ to be the linear span of
\begin{eqnarray}
\bar{f}_{g_{0}}(u)-u\;\;\;\mbox{ for }g_{0}\in G_{0},\; u\in U.
\end{eqnarray}

\bl{lideal}
The defined subspace $I$ is an ideal of $U$ with
$I\cap V^{0}=0$. Furthermore, $I=0$ if and only if $G_{0}=0$.
\el

\pf Let $g_{0}\in G_{0},\; h,h'\in G$ and let $u\in V^{h},\; v\in V^{h'}$.
Then
\begin{eqnarray}\label{eabstract1}
\bar{f}_{g_{0}}(u)=\Res_{z_{2}}\bar{f}_{g_{0}}(Y(u,z_{2}){\bf 1})
=\Res_{z_{2}}Y(u,z_{2})\bar{f}_{g_{0}}({\bf 1}).
\end{eqnarray}
Since
\begin{eqnarray}\label{eabstract2}
Y(v,z)\bar{f}_{g_{0}}({\bf 1})=\bar{f}_{g_{0}}(Y(v,z){\bf 1})\in U[[z]],
\end{eqnarray}
we have (cf. [Li2])
\begin{eqnarray}\label{eabstract3}
Y(v,z_{0}+z_{2})Y(u,z_{2})f_{g_{0},h}({\bf 1})=
Y(Y(v,z_{0})u,z_{2})f_{g_{0},h}({\bf 1}).
\end{eqnarray}
Using (\ref{eabstract1})-(\ref{eabstract3}) we obtain
\begin{eqnarray}
& &Y(v,z_{0})(\bar{f}_{g_{0}}(u)-u)\nonumber\\
&=&\Res_{z_{2}}Y(v,z_{0}+z_{2})Y(u,z_{2})f_{g_{0},h}({\bf 1})-Y(v,z_{0})u
\nonumber\\
&=&\Res_{z_{2}}Y(Y(v,z_{0})u,z_{2})f_{g_{0},h}({\bf 1})-Y(v_{0},z)u\nonumber\\
&=&f_{g_{0},h}(Y(v,z_{0})u)-Y(v,z_{0})u\nonumber\\
&=&\bar{f}_{g_{0}}(Y(v,z_{0})u)-Y(v,z_{0})u.
\end{eqnarray}
It follows immediately that $I$ is an ideal.

Clearly, $I\ne 0$ if $G_{0}\ne \{0\}$, so $I=0$ implies $G_{0}=0$.
If $G_{0}=0$, with $f_{0}=1$ from the definition of $I$ we have $I=0$.
$\;\;\;\;\Box$

\bp{pabstract2}
The algebra $U$ is simple if and only if $G_{0}=0$. Furthermore,
the quotient algebra $\bar{U}=U/I$ is simple. 
\ep

\pf From Lemma \ref{lideal}, $U$ is not simple if $G_{0}\ne 0$.
Now we prove that $U$ is simple if $G_{0}=\{0\}$.
In view of Lemma \ref{labstract1}, all $V^{g}$ for $g\in G$ are
non-isomorphic irreducible $V^{0}$-modules. Then any nonzero ideal
of $U$ as a $V^{0}$-module must be a sum of some $V^{g}$.
Then it follows from the conditions (1)-(3) that
any nonzero ideal of $U$ must be $U$. That is, $U$ is simple.

Note that $\bar{U}=U/I$ is a vertex operator (super)algebra
graded by group $G/G_{0}$ with all the conditions (1)-(3) being
satisfied. Furthermore, $\bar{U}$ is a direct sum of non-isomorphic
simple current $V^{0}$-modules. From Part one, $\bar{U}$ must be simple.
$\;\;\;\;\Box$

Applying Proposition \ref{pabstract2} to $V[L]$
we immediately have (cf. [DLM2], Corollary 3.13):

\bc{cdlm2}
Let $V, L$ be as before. Set
\begin{eqnarray}
L_{0}=\{\a\in L\;|\; V^{(\a)}\simeq V\;\;\;\mbox{ as } V\mbox{-modules}\}.
\end{eqnarray}
Then $V[L]$ has an ideal $I$ such that $V[L]/I$ is simple with $V$ as 
a subalgebra and such that as a $V$-module
\begin{eqnarray}
V/I\simeq \oplus_{\a\in S}V^{(\a)},
\end{eqnarray}
where $S$ is a complete set of representatives of cosets of $L_{0}$ in $L$.
$\;\;\;\;\Box$
\ec

\subsection{Multiplicity-free extension $V[L]$}
In this subsection we consider extended vertex (super)algebra 
$V[L]$ in which each $V$-module $V^{(\a)}$
is multiplicity-free. We shall classify all irreducible 
$V[L]$-modules in terms of irreducible $V$-modules
and determine the fusion rules of $V[L]$-modules 
by the fusion rules of $V$-modules. Our concrete 
examples we shall construct in Section 4 is of this type,
so that all the results of this subsection apply to
those examples.

Throughout this subsection we assume that for any irreducible 
$V$-module $W$ and for $\a,\beta\in L$,
$W^{(\a)}\simeq W^{(\beta)}$ as $V$-modules
if and only if $\a=\beta$.

We shall need the following result:

\bp{puniqueness}
The vertex (super)algebra $V[L]$ is simple. Furthermore,
if $Y_{1}$ and $Y_{2}$ are
two simple vertex operator (super)algebra 
structures on $V[L]$ extending the $V$-module structure $Y_{V}$,
then vertex operator (super)algebras $(V[L],Y_{1})$ and 
$(V[L],Y_{2})$ are isomorphic.	
\ep

\pf First, we prove that $V[L]$ is simple.
Notice that in the definition (\ref{edefvoa}) of the vertex operator map $Y$,
$\psi_{\a+\beta}, \psi_{-\a}, \psi_{-\beta},
E^{-}(\alpha,z)$, $\Delta(\beta,z)$ and $\Delta(\alpha,-z)$ 
are invertible elements and that
$$Y(a,z)b\in V^{(\a+\beta)}\;\;\;\mbox{ for }a\in V^{(\a)},\;b\in V^{(\beta)}$$
(cf. (\ref{edefvoa})).
Since $V$ is simple, $Y(u,z)v\ne 0$ for $0\ne u, v\in V$ 
([DL], Proposition 11.9, or [FHL], Remark 5.4.6).
Then it follows from Proposition \ref{pabstract2} 
immediately that $V[L]$ is simple.

For $\a,\beta\in L$, because $V^{(\a)}$ and $V^{(\beta)}$ are simple currents,
there exists $\epsilon'(\a,\beta)\in \C^{\times}$ such that
\begin{eqnarray}
Y_{2}(a,z)b=\epsilon'(\a,\beta)Y_{1}(a,z)b\;\;\;\mbox{ for }
a\in V^{(\a)},\; b\in V^{(\beta)}.
\end{eqnarray}
It follows from weak associativity of vertex operators 
(cf. (\ref{e3.54})) that
\begin{eqnarray}
\epsilon'(\a,\beta+\gamma)\epsilon'(\beta,\gamma)
=\epsilon'(\a,\beta)\epsilon'(\a+\beta,\gamma)
\end{eqnarray}
for $\a,\beta\in L$. That is, $\epsilon'$ is a ($\C^{\times}$-valued)
2-cocycle on $L$.
We also have
\begin{eqnarray}
\epsilon'(0,\a)=\epsilon'(\a,0)=1.
\end{eqnarray}
Since $V=V^{(0)}$ is even for both superalgebra structures,
each structure corresponds a sublattice $L_{i}$ of $L$ of index 2
with $V[L_{i}]$ being the even parts.

Now we claim that $L_{1}=L_{2}$.
Otherwise, suppose $L_{1}-L_{2}\ne \emptyset$ and let
$\beta\in L_{1}-L_{2}$. Then we have the following skew-symmetry 
\begin{eqnarray}
& &Y_{1}(a,z)b=e^{zL(-1)}Y_{1}(b,-z)a,\\
& &Y_{2}(a,z)b=-e^{zL(-1)}Y_{2}(b,-z)a
\end{eqnarray}
for $a,b\in V^{(\b)}$. Since both $Y_{1}$ and $Y_{2}$ extend 
the $V$-module structure, the two vertex superalgebra structures have
the same Virasoro vector. Consequently,
\begin{eqnarray}
\epsilon'(\b,\b)=-\epsilon'(\b,\b),
\end{eqnarray}
which is impossible because $\epsilon'(\b,\b)\ne 0$.

With $L_{1}=L_{2}$, using the skew-symmetry we obtain
\begin{eqnarray}
\epsilon'(\a,\beta)=\epsilon'(\beta,\a)\;\;\;\mbox{ for }\a,\beta\in L.
\end{eqnarray}
It follows from the proof of Propositions 5.1.2 and 5.2.3 in [FLM]
(with $\Z/s\Z$ being replaced by $\C^{\times}$) that
$\epsilon'$ is a $2$-coboundary. Then the two vertex superalgebra 
structures on $V[L]$  are equivalent.
$\;\;\;\;\Box$

\br{runiquenessalgebra}
{\em More generally, let $G$ be a (finite or infinite) abelian group
and let $V$ be a simple vertex operator algebra and
$V[G]=\oplus_{g\in G}V^{g}$ be a $V$-module with each $V^{g}$
 being an irreducible $V$-submodule. Furthermore, assume that 
each $V^{g}$ is a simple current of $V$. Then the set of equivalence class of
simple vertex operator (super)algebra structures on $V[G]$ which
extend the $V$-module one-to-one corresponds to the set of 
equivalence classes of
symmetric $\C^{\times}$-valued $2$-cocycles of $G$.}
\er

Similar to Proposition \ref{puniqueness} we have (cf. [DLM2], Lemma 4.2):

\bp{pmoduleirreducibility}
Let $W$ be an irreducible $V$-module.
Then $W[L]$ is irreducible. $\;\;\;\;\;\Box$
\ep

The following theorem gives the complete reducibility
for every $V[L]$-module under certain conditions:

\bt{textensionregularity}
Assume that there is a sublattice $L_{1}$ of $L$ such that $V[L_{1}]$
is regular and that every irreducible $V[L_{1}]$-module is a direct sum of
irreducible $V$-modules. Let $\sigma$ be an automorphism of $V[L]$
such that $\sigma$ fixes $V[L_{1}]$ point-wise.
Then any $\sigma$-twisted weak $V[L]$-module is a direct sum of
irreducible $\sigma$-twisted $V[L]$-modules of type $W[L]$ with 
$\sigma=\sigma_{W}$.
In particular, $V[L]$ is regular.
\et

\pf Let $M$ be a $\sigma$-twisted weak $V[L]$-module.
Since $\sigma$ fixes $V[L_{1}]$ point-wise,
$M$ is a weak $V[L_{1}]$-module. With $V[L_{1}]$ being regular,
$M$ is a direct sum of irreducible (ordinary) 
$V[L_{1}]$-modules. With the assumption that
each irreducible $V[L_{1}]$-module
is a direct sum of irreducible $V$-modules,
$M$ is a direct sum of irreducible $V$-modules.
Let $W$ be an irreducible $V$-submodule of $M$.
By Proposition \ref{pextensionhom2}, $\sigma=\sigma_{W}$ and 
there exists a $V[L]$-homomorphism $\phi[L]$ from $W[L]$ to $M$
extending the embedding $\phi$ of $W$ into $M$.
Since $W[L]$ is $V[L]$-irreducible 
(Proposition \ref{pmoduleirreducibility}), 
the $V[L]$-submodule of $M$
generated by $W$, which is the image of $\phi[L]$, is an 
irreducible $\sigma$-twisted $V[L]$-module.
Therefore, $M$ is a direct sum of
irreducible $\sigma$-twisted $V[L]$-modules of type $W[L]$. $\;\;\;\;\Box$

Let $W_{1},W_{2},W_{3}$ be irreducible $V$-modules such that
$\sigma_{W_{i}}=1$, or equivalently, the extensions $W_{1}[L], W_{2}[L]$ 
and $W_{3}[L]$ are $V[L]$-modules. 
For each $\a\in L$, with $W_{3}^{(\a)}$ being a $V$-submodule of 
$W_{3}[L]$, by Proposition \ref{pextensionhom2},
there is a $V[L]$-homomorphism $g_{\a}$ from $W_{3}^{(\a)}[L]$ 
to $W_{3}[L]$.
Then with Proposition \ref{pextensionintertwining}
 we obtain a natural linear map $f_{\a}$ from 
$I_{W_{1}W_{2}}^{W_{3}^{(\a)}}$ to $I_{W_{1}[L]W_{2}[L]}^{W_{3}[L]}$
defined by
\begin{eqnarray}
f_{\a}(I)=g_{\a}\circ I[L].
\end{eqnarray}
The next result gives a precise connection between 
the fusion rules for $V$-modules
and the fusion rules for $V[L]$-modules:

\bt{tfusionrule}
Let $W_{1},W_{2},W_{3}$ be irreducible $V$-modules such that
$\sigma_{W_{i}}=1$, which implies that
$W_{1}[L], W_{2}[L], W_{3}[L]$ are irreducible $V[L]$-modules.
In addition, we assume that $V$ is quasi-rational.
Then $f=\coprod_{\a\in L}f_{\a}$ is a linear isomorphism from
$\coprod_{\a\in L}I_{W_{1}W_{2}}^{W_{3}^{(\a)}}$ to
$I_{W_{1}[L]W_{2}[L]}^{W_{3}[L]}$.
In particular,
\begin{eqnarray}
N_{W_{1}[L]W_{2}[L]}^{W_{3}[L]}=\sum_{\a\in L}N_{W_{1}W_{2}}^{W_{3}^{(\a)}}.
\end{eqnarray}
\et

\pf Since $W_{3}^{(\a)}\simeq W_{3}^{(\b)}$ only if $\a=\b$,
clearly, $f_{\a}$ is one-to-one. On the other hand,
if $\cal{Y}$ is an intertwining operator of type
${W_{3}[L]\choose W_{1}[L]W_{2}[L]}$. Then by restricting $\cal{Y}$ to
$W_{1}\times W_{2}$ we have an intertwining operator
$I$ of type ${W_{3}^{(\a)}\choose W_{1}W_{2}}$ for a unique $\a\in L$.
It is clear that ${\cal{Y}}=I[L]$. This completes the proof.
$\;\;\;\;\Box$

We now describe the Verlinde algebra of $V[L]$ in terms of
the Verlinde algebra of $V$ explicitly.
Let ${\cal{A}}(V)$ be the Verlinde algebra of $V$.
The Verlinde algebra ${\cal{A}}(V[L])$ with a basis
$[W[L]]$ for $[W]\in {\cal{A}}(V)$ with $\sigma_{W}=1$.
First, we have:

\bl{lverlindesubalgebra}
All $[W]$ with $\sigma_{W}=1$ linearly span
a subalgebra $A(V,L)$ of ${\cal{A}}(V)$.
\el

\pf Suppose that $[W_{1}], [W_{2}]\in \Irr(V)$ with 
$\sigma_{W_{1}}=\sigma_{W_{2}}=1$.
Let $h_{1}\in P(W_{1}),\; h_{2}\in P(W_{2})$.
Then $\sigma_{W_{1}}=\sigma_{W_{2}}=1$ amount to $h_{1},h_{2}\in L^{o}$.
Let ${\cal{Y}}$ be a nonzero intertwining operator
of type ${W_{3}\choose W_{1}W_{2}}$ for some irreducible 
$V$-module $W_{3}$. Then
$h_{1}+h_{2}\in P(W_{3})$ because for $h\in H,\; w_{(i)}\in W_{i}^{(0,h_{i})}$,
\begin{eqnarray}
h(0){\cal{Y}}(w_{(1)},z)w_{(2)}&=&{\cal{Y}}(h(0)w_{(1)},z)w_{(2)}+
{\cal{Y}}(w_{(1)},z)h(0)w_{(2)}\nonumber\\
&=&B(h,h_{1}+h_{2}){\cal{Y}}(w_{(1)},z)w_{(2)}.
\end{eqnarray}
Since $h_{1}+h_{2}\in L^{o}$, $\sigma_{W_{3}}=1$.
Then $[W_{3}]\in A(V,L)$. The proof is complete.
$\;\;\;\;\Box$

Define a subspace $R$ of $A(V,L)$ linearly spanned by
\begin{eqnarray}
[W]-[W^{(\a)}]\;\;\;\;\;\mbox{ for }\a\in L.
\end{eqnarray}
Then $R$ is an two-sided ideal of $A(V,L)$.
Indeed, let $W_{1}$ and $W_{2}$ be irreducible $V$-modules with 
$\sigma_{W_{i}}=1$ for $i=1,2$. For $\a\in L$, by Proposition 2.10,
\begin{eqnarray}
I_{W_{1}W}^{W_{2}}\simeq I_{W_{1}W^{(\a)}}^{W_{2}^{(\a)}}.
\end{eqnarray}
Thus
\begin{eqnarray}
[W_{1}]\cdot ([W]-[W^{(\a)}])=\sum_{[W_{2}]\in Irr(V)}N_{W_{1}W}^{W_{2}}
([W_{2}]-[W_{2}^{(\a)}])\in R.
\end{eqnarray}
Since ${\cal{A}}(V)$ is a commutative algebra, $R$ is a two-sided ideal.
Furthermore, by Proposition \ref{pequivalencew[l]},
 $\overline{[W_{1}]}=\overline{[W_{2}]}$ in the quotient algebra
$A(V,L)/R$ if and only if
$[W_{1}[L]]=[W_{2}[L]]$ in ${\cal{A}}(V[L])$.
Then in view of Theorem  \ref{tfusionrule} we immediately have:

\bp{pverlinealgebra}
The subspace $R$ is a two-sided ideal of $A[V,L]$
and the Verlinde algebra ${\cal{A}}(V[L])$ is canonically
isomorphic to the quotient algebra of $A(V,L)$ modulo $R$.$\;\;\;\;\Box$
\ep

\section{Extended vertex operator (super)algebras of affine types}
In this section we shall specialize the vertex (super)algebra $V[L]$
constructed in Section 3 from a pair $(V,L)$
to obtain extensions of vertex operator algebras
associated with affine Lie algebras $\hat{\fg}$. 
In the case $\fg=sl(2)$, we shall obtain Feigin-Miwa's 
extended vertex operator (super)algebras $A_{k}$ [FM].
To apply the results of Section 3 we need to define $V$ and $L$ explicitly
and check the necessary conditions. 
For each type, $V$ will be the tensor product vertex operator algebra
$L_{\fg }(k,0)\otimes M_{{\bf h}'}(1,0)$ where ${\bf h}'$ is 
a 1 or 2-dimensional vector space equipped with a
nondegenerate symmetric bilinear form.
When defining
${\bf h}'$, we follow two basic principles: (1) To include
all the $L(k,0)$-simple currents in the construction of $V[L]$. (2) To make
$\dim {\bf h}'$, or equivalently, to make the rank of $V[L]$
as small as possible. 
After ${\bf h}'$ is chosen, we still have
plenty of choices for the bilinear form on ${\bf h}'$.
Another principle to follow is to
make $V[L]$ as large as possible. 
We shall define ${\bf h}'$ and $L$ type by type. 

\subsection{A complete reducibility theorem for a certain family of $V[L]$}

In Section 3, for a general pair $(V,L)$, under a certain assumption we 
proved a complete reducibility theorem (Theorem \ref{textensionregularity})
for extended algebra $V[L]$.
In this section, we shall consider a family of $V[L]$
such that the assumption of Theorem \ref{textensionregularity} holds.
All the extended algebras we shall construct later
belong to this family, so the complete reducibility theorem holds 
for all of them.

Let $\fg=\fg_{1}\oplus\cdots\oplus\fg_{r}$ be a semisimple Lie algebra 
with a Cartan subalgebra
\begin{eqnarray}
{\bf h}={\bf h}_{1}+\cdots +{\bf h}_{r},
\end{eqnarray}
where $\fg_{i}$ are simple factors with Cartan subalgebras ${\bf h}_{i}$,
equipped with the normalized Killing forms. 
Let $Q^{\vee}=Q_{1}^{\vee}+\cdots +Q_{r}^{\vee}$ 
and $P^{\vee}=P^{\vee}_{1}+\cdots +P_{r}^{\vee}$ be the coroot lattice and
coweight lattice of $\fg$, respectively, where $Q_{i}^{\vee}$ and 
$P_{i}^{\vee}$ are the coroot lattice and
coweight lattice of $\fg_{i}$.

Let ${\bf k}=(k_{1},\dots,k_{r})$ be an $r$-tuple of nonnegative integers.
Set
\begin{eqnarray}
L_{\fg}({\bf k},0)
=L_{\fg_{1}}(k_{1},0)\otimes \cdots\otimes L_{\fg_{r}}(k_{r},0),
\end{eqnarray}
equipped with the standard tensor product vertex operator algebra 
structure. Set
\begin{eqnarray}
P_{{\bf k}}=\{(\lambda^{1},\dots,\lambda^{r})\;|\; 
\lambda^{i}\in P_{k_{i}}(\fg_{i})\},
\end{eqnarray}
where $P_{k_{i}}(\fg_{i})$ stands for $P_{k_{i}}$ for 
the Lie algebra $\fg_{i}$.
Then $L_{\fg}({\bf k},0)$ is regular [DLM1], i.e., every weak module is
a direct sum of irreducible (ordinary) modules 
$L_{\fg}({\bf k},{\bf \lambda})$, where
\begin{eqnarray}
L_{\fg}({\bf k},{\bf \lambda})
=L_{\fg_{1}}(k_{1},\lambda^{1})\otimes \cdots\otimes 
L_{\fg_{r}}(k_{r},\lambda^{r})
\end{eqnarray} 
for ${\bf \lambda}=(\lambda^{1},\dots,\lambda^{r})\in P_{\bf k}$.

Let ${\bf h}'$ be a finite-dimensional vector space equipped with a
nondegenerate symmetric bilinear form $\<\cdot,\cdot\>$.
Associated to ${\bf h}'$ is the vertex operator algebra
$M_{{\bf h}'}(1,0)$. Set
\begin{eqnarray}\label{e4.4}
V=L_{\fg}({\bf k},0)\otimes M_{{\bf h}'}(1,0),
\end{eqnarray}
equipped with the standard tensor product vertex operator algebra 
structure. The algebra $V$ can be considered as the vertex operator algebra
associated to the affine algebra of the reductive Lie algebra $\fg+{\bf h}'$.

Set
\begin{eqnarray}\label{e4.5}
H={\bf h}+{\bf h}'\subset \fg+{\bf h}'=V_{(1)}.
\end{eqnarray}
Then $H$ is a Cartan subalgebra of the reductive Lie algebra $\fg+{\bf h}'$.
Clearly, (\ref{e3.1}) holds and
\begin{eqnarray}
B(h_{1},h_{2})=\delta_{i,j}k_{i}\<h_{1},h_{2}\> \;\;\;
\mbox{ for }h_{1}\in {\bf h}_{i},\;h_{2}\in {\bf h}_{j}
\end{eqnarray}
because $h_{1}(1)h_{2}=h_{1}(1)h_{2}(-1){\bf 1}
=\delta_{i,j}k_{i}\<h_{1},h_{2}\>{\bf 1}$
(and $h_{1}(1)h_{2}=B(h_{1},h_{2}){\bf 1}$). We also have
\begin{eqnarray}
B(h_{1}',h_{2}')=\<h_{1}',h_{2}'\> \;\;\;
\mbox{ for }h_{1}', h_{2}'\in {\bf h}'.
\end{eqnarray}
For $h\in H$, we write
\begin{eqnarray}
h=h^{1}+\cdots +h^{r}+h',\;\;\;\; h=h''+h',
\end{eqnarray}
where $h^{i}\in {\bf h}_{i},\; h'\in {\bf h}',\;h''\in {\bf h}$.

For ${\bf \lambda}=(\lambda^{1},\dots,\lambda^{r}),\; \gamma\in {\bf h}'$
with $\lambda^{i}\in P_{k_{i}}(\fg_{i})$, we set
\begin{eqnarray}
W({\bf \lambda},\gamma)
=L_{\fg_{1}}(k_{1},\lambda^{1})\otimes \cdots \otimes 
L_{\fg_{r}}(k_{r},\lambda^{r})\otimes M_{{\bf h}'}(1,\gamma).
\end{eqnarray}
Then from [FHL] all such $W({\bf \lambda},\gamma)$ form a complete set 
of representatives of equivalence classes of irreducible $V$-modules.

Let $L$ be a subgroup of $H$ such that $(L,B)$ is an integral
lattice of finite rank. 
To have a vertex (super)algebra $V[L]$, 
we shall also need the condition that $\a(0)$ has only integral 
eigenvalues on $V$ for every $\a\in L$. 
For $\a\in L$, since $\a'(0)$ 
acts as zero on $V$, $\a(0)$ has only integral eigenvalues on $V$
if and only if $\a''(0)$ has only integral eigenvalues on $L_{\fg}({\bf k},0)$.
Then we immediately have:

\bl{lintegralvalues} 
For $\a\in L$, $\a(0)$ has only integral eigenvalues on $V$ if and only if
$\a''\in P^{\vee}$, the coweight lattice of $\fg$.$\;\;\;\;\Box$
\el

Set
\begin{eqnarray}
L'=\{\a'\;|\; \a\in L\},\;\;\;L''=\{\a''\;|\; \a\in L\}.
\end{eqnarray}

\bl{lfamily}
Assume that the projection of $L$ onto $L'$ is one-to-one.
Then for ${\bf \lambda}\in P_{\bf k},\; \gamma\in {\bf h}'$ and
for $\a,\b\in L$,
$W({\bf \lambda},\gamma)^{(\a)}\simeq W({\bf \lambda},\gamma)^{(\b)}$
if and only if $\a=\b$.
\el

\pf Because $\Delta(\a,z)=\Delta(\a'',z)\Delta(\a',z)$ and
$\Delta(\a'',z)=1$ on $M_{{\bf h}'}(1,\gamma)$ and $\Delta(\a',z)=1$ on 
$L_{\fg}({\bf k},{\bf \lambda})$ for $\a\in L$, we have
\begin{eqnarray}
W({\bf \lambda},\gamma)^{(\a)}
=L_{\fg}({\bf k},{\bf \lambda})^{(\a'')}\otimes M_{{\bf h}'}(1,\gamma)^{(\a')}
\simeq 
L_{\fg}({\bf k},{\bf \lambda})^{(\a'')}\otimes M_{{\bf h}'}(1,\gamma+\a').
\end{eqnarray}
We knew $M_{{\bf h}'}(1,\gamma+\a')\simeq M_{{\bf h}'}(1,\gamma+\b')$
if and only if $\a'=\b'$. Then it follows immediately. $\;\;\;\;\Box$

Now we have:

\bp{pgeneralaffine}
Let $V, H$ be defined as in (\ref{e4.4}) and (\ref{e4.5})
and let $L$ be a subgroup of $H$ such that $(L,B)$ is an integral 
lattice of finite rank. Assume that $L''\subset P^{\vee}$, 
the projection of $L$ onto $L'$ is one-to-one and 
$L'$ is a positive definite lattice.
Then $V[L]$ equipped with the vertex operator 
map $Y$ defined in (\ref{edefvoa})
is a simple vertex operator (super)algebra.
\ep

\pf Since $L''\subset P^{\vee}$, by Lemma \ref{lintegralvalues},
for $\a\in L$, $\a(0)$ acting on $V$ has only 
integral eigenvalues. With Lemma \ref{lfamily}, 
in view of Corollary \ref{c3.13} and Proposition \ref{puniqueness},
$V[L]$ is a simple vertex (super)algebra with all $L(0)$-weights being 
half integers.
Now we only need to verify the two grading restrictions [FLM].

For $\a\in L$, we have
\begin{eqnarray}
V^{(\a)}=L_{\fg}({\bf k},0)^{(\a'')}\otimes M_{{\bf h}'}(1,0)^{(\a')}
=L_{\fg}({\bf k},0)^{(\a'')}\otimes M_{{\bf h}'}(1,\a').
\end{eqnarray}
With $L_{\fg}({\bf k},0)$ being regular, $L_{\fg}({\bf k},0)^{(\a'')}$ is
an irreducible $L_{\fg}({\bf k},0)$-module, which  is unitary.
Then the $L(0)$-weights of
$L_{\fg}({\bf k},0)^{(\a'')}$ are nonnegative. On the other hand,
the lowest weight of $M_{{\bf h}'}(1,\a')$ is ${1\over 2}\<\a',\a'\>$.
Thus each $V^{(\a)}$
satisfies the two grading restrictions with the lowest weight
at least ${1\over 2}\<\a',\a'\>$.
Since
the projection of $L$ onto $L'$ is one-to-one and $L'$ is positive definite,
for every $n\in {1\over 2}\Z$, $V^{(\a)}_{(n)}\ne 0$ only for
finitely many $\a$.
Then the two grading restrictions follows immediately.
$\;\;\;\;\Box$

Furthermore, we have:

\bt{tregularitysemisimple}
Let $V, L$ be as in Proposition \ref{pgeneralaffine} with all 
the assumptions. In addition we assume that $\dim {\bf h}'=\mbox{ rank}\; L'$.
Let $\sigma$ be an automorphism of $V[L]$ of finite order
which fixes $V$ point-wise.
Then every $\sigma$-twisted weak $V[L]$-module is a direct sum
of irreducible (ordinary) $\sigma$-twisted $V[L]$-modules
isomorphic to $W({\bf \lambda},\gamma)[L]$ with 
$\sigma=\sigma_{W({\bf \lambda},\gamma)}$.
In particular, every weak $V[L]$-module
is a direct sum of irreducible (ordinary) $V[L]$-modules
isomorphic to $W({\bf \lambda},\gamma)[L]$ for ${\bf \lambda}\in P_{\bf k}$,
$\gamma\in (L')^{o}$ (the dual of $L'$) with 
\begin{eqnarray}\label{e4.13}
\lambda^{1}(\a^{1})+\cdots+\lambda^{r}(\a^{r})+\<\gamma,\a'\>\in \Z
\;\;\;\mbox{ for }\a\in L.
\end{eqnarray}
\et

\pf Denote by $o(\sigma)$ the order of $\sigma$. 
In view of Lemma \ref{lfamily},
$V^{(\a)}$ for $\a\in L$ are non-isomorphic irreducible $V$-modules.
Then $\sigma(V^{(\a)})=V^{(\a)}$ for $\a\in L$ and
$\sigma$ acts on $V^{(\a)}$ as a scalar, 
which is an $o(\sigma)$-th root of unity. 
Therefore, $\sigma$ acts trivially on $V^{(mo(\sigma)\a)}$ 
for $\a\in L,\; m\in \Z$.

Because $L_{\fg}({\bf k},0)$ has only finitely many irreducible modules 
up to equivalence, there exists a positive integer $d_{1}$ such that
as $L_{\fg}({\bf k},0)$-modules,
\begin{eqnarray}
L_{\fg}({\bf k},0)^{(d_{1}\a'')}\simeq L_{\fg}({\bf k},0)\;\;\;\mbox{ for all }
\a\in L.
\end{eqnarray}
Let $d_{2}$ be another positive integer such that $d_{2}L'$ is 
an even lattice. Set $d=o(\sigma)d_{1}d_{2}$ and $L_{1}=dL$. Then
\begin{eqnarray}
V[L_{1}]=\oplus_{\a\in L}V^{(d\a)}.
\end{eqnarray}
Then $V[L_{1}]$ is a simple vertex operator subalgebra of $V[L]$ 
and $\sigma$ fixes 
$V[L_{1}]$ point-wise.
Furthermore, as a $V$-module,
\begin{eqnarray}
V^{(d\a)}
\simeq L_{\fg}({\bf k},0)^{(d\a'')}\otimes M_{{\bf h}'}(1,d\a')
\simeq L_{\fg}({\bf k},0)\otimes M_{{\bf h}'}(1,d\a')
\end{eqnarray}
for $\a\in L$, hence
\begin{eqnarray}
V[L_{1}]\simeq L_{\fg}({\bf k},0)\otimes V_{dL'}.
\end{eqnarray}
(Here we used the fact that $\dim {\bf h}'=\mbox{rank}\; L'$.)
Note that $L_{\fg}({\bf k},0)\otimes V_{dL'}$ is a natural simple 
vertex operator algebra, which is regular.
It follows from Proposition \ref{puniqueness} that
$V[L_{1}]$ is regular. Clearly, each irreducible $V[L_{1}]$-module
is a direct sum of irreducible $V$-modules.
Then it follows from Theorem \ref{textensionregularity}
immediately that every $\sigma$-twisted weak $V[L]$-module is a direct sum of 
irreducible (ordinary) $\sigma$-twisted $V[L]$-modules of type 
$W({\bf \lambda},\gamma)$
with $\sigma=\sigma_{W({\bf \lambda},\gamma)}$.

{}From Lemma \ref{lautomorphism}, $\sigma_{W({\bf \lambda},\gamma)}=1$ 
if and only if for $\a\in L$,
$\a (0)$ has only integral eigenvalues on  $W({\bf \lambda},\gamma)$.
With $(\lambda^{1},\dots,\lambda^{r},\gamma)$ being an $H$-weight of
$W({\bf \lambda},\gamma)$, we see that $\sigma_{W({\bf \lambda},\gamma)}=1$ 
if and only if 
\begin{eqnarray}
\lambda^{1}(\a^{1})+\cdots +\lambda^{r}(\a^{r})+\<\gamma,\a'\>\in \Z
\;\;\;\mbox{ for }\a\in L,
\end{eqnarray}
which furthermore implies that $\gamma\in (L')^{o}$ because
${\bf \lambda}\in P_{\bf k}$.
This completes the proof.$\;\;\;\;\Box$

\br{rgeneralU}
{\em From the proof of Theorem \ref{tregularitysemisimple}, 
one can easily see that
the regularity result still holds for $V[L]$ if we replace
$L_{\fg}({\bf k},0)$ by any regular vertex operator algebra $U$.}
\er

Because 
$$W({\bf \lambda},\gamma)^{(\a)}
=L({\bf k},{\bf \lambda})^{(\a'')}\otimes M_{{\bf h}'}(1,\gamma+\a')$$
and $M_{{\bf h}'}(1,\gamma)\simeq M_{{\bf h}'}(1,\gamma')$
if and only if $\gamma=\gamma'$,
in view of Proposition \ref{pequivalencew[l]}, 
we see that 
$W({\bf \lambda},{\bf \gamma})[L]\simeq W({\bf \lambda}',{\bf \gamma}')[L]$
if and only if there is $\a\in L$ such that
\begin{eqnarray}
\gamma'=\gamma+\a',\;\;\;
L({\bf k},{\bf \lambda}')^{(\a'')}\simeq L({\bf k},{\bf \lambda}).
\end{eqnarray}
To describe explicitly the equivalence relation on the set of 
$W({\bf \lambda},{\bf \gamma})[L]$, or to get
a complete set of equivalence
classes of irreducible $V[L]$-modules, we need to know
$L({\bf k},{\bf \lambda}')^{(\a'')}$ as a $\hat{\fg}$-module.
Of course, from Theorem \ref{tli3},
\begin{eqnarray}
[L({\bf k},{\bf \lambda}')^{(\a'')}]
=[L({\bf k},{\bf \lambda}')]\cdot [L({\bf k},0)^{(\a'')}].
\end{eqnarray}

Nevertheless, in view of Proposition \ref{pverlinealgebra} we immediately have:

\bp{pverlindegeneral}
The subspace $A$ of ${\cal{V}}(L({\bf k},0))\otimes \C[(L')^{o}]$,
linearly spanned by
\begin{eqnarray}
[L({\bf k},{\bf \lambda})]\otimes e^{\gamma}
\end{eqnarray}
for ${\bf \lambda}\in P_{\bf k},\; \gamma\in (L')^{o}$ satisfying
(\ref{e4.13}), is a subalgebra. Furthermore,
the Verlinde algebra ${\cal{V}}(V[L])$
is canonically isomorphic to the 
quotient algebra of $A$ modulo the relations:
\begin{eqnarray}
[L({\bf k},{\bf \lambda})]\otimes e^{\gamma}
-[L({\bf k},{\bf \lambda})^{(\a'')}]\otimes e^{\gamma+\a'}
\end{eqnarray}
for $\a\in L$.$\;\;\;\;\Box$
\ep

\subsection{Extended vertex operator (super)algebras $A_{k}$ of type $sl(n+1)$}
Starting from this subsection we shall work on the 
setting  of Section 4.1 and we shall consider 
a simple Lie algebra $\fg$, i.e., $r=1$.
For $\fg$ of each type, we take
\begin{eqnarray}\label{edvsln+1}
V=L(k,0)\otimes M_{{\bf h}'}(1,0)
\end{eqnarray}
and we define $A_{k}(\fg)$ to be the extended algebra $V[L]$
for a certain $L$.
We shall case by case define the pair $({\bf h}', \<\cdot,\cdot\>)$ 
and the lattice $L$, and then verify that
$(L,B)$ is an integral lattice, $L''\subset P^{\vee}$,
$L'$ is positive-definite and the projection of $L$ onto $L'$ is
one-to-one, so that Proposition \ref{pgeneralaffine} and 
Theorem \ref{tregularitysemisimple} hold.

In this subsection we shall consider the case $\fg=sl(n+1)$. 
For a fixed positive integer $k$, $L(k,0)$
has $n+1$ simple currents $L(k,k\lambda_{i})$ for $i=0,\dots,n$.
By Corollary \ref{csimplecurrentfusion} 
the equivalence classes of the $(n+1)$ simple currents
form a cyclic group of order $(n+1)$ with 
$[L(k,k\lambda_{1})]$ as a generator.

Recall that $h^{(i)}\in {\bf h}$ with $\a_{j}(h^{(i)})=\delta_{i,1}$
for $i,j=1,\dots, n$. From [H] (Table 1 on page 69) and 
Lemma \ref{lsimplefact}, we have
\begin{eqnarray}
& &h^{(i)}={1\over n+1}\left((n+1-i)\a_{1}^{\vee}+2(n+1-i)\a_{2}^{\vee}+\cdots 
+(i-1)(n+1-i)\a_{i-1}^{\vee}\right)\nonumber\\
& &\;\;\;\;\;\;+ {1\over n+1}\left(i(n+1-i)\a_{i}^{\vee} 
+i(n-i)\a_{i+1}^{\vee}+
\cdots +i\a_{n}^{\vee}\right),\label{ehi}\\
& &\<h^{(i)},h^{(i)}\>={i(n+1-i)\over n+1}.
\end{eqnarray}

Define ${\bf h}'=\C\alpha'$ 
to be a one-dimensional vector space equipped with a symmetric bilinear 
form $\<\cdot,\cdot\>$ such that
\begin{eqnarray}
\<\a',\a'\>=\frac{k}{n+1}.
\end{eqnarray}
Set
\begin{eqnarray}\label{edlsln+1}
L=\Z \a,\;\;\;\mbox{ where }\a=h^{(1)}+\alpha'.
\end{eqnarray}
Because
\begin{eqnarray}
B(\a,\a)=B(h^{(1)},h^{(1)})+B(\a',\a')={kn\over n+1}+{k\over n+1}=k\in \Z,
\end{eqnarray}
$(L,B)$ is an integral lattice. Clearly, $L''=\Z h^{(1)}\subset P^{\vee}$,
$L'$ is positive-definite and the projection of $L$ onto $L'$ is one-to-one.
By Proposition \ref{pgeneralaffine},
we have a simple vertex operator (super)algebra $V[L]$.

\bd{dsln+1}
{\em We define $A_{k}(sl(n+1))$ to be the simple vertex operator 
(super)algebra $V[L]$
with $V$ and $L$ being defined in (\ref{edvsln+1}) and (\ref{edlsln+1}).}
\ed

\br{rother}
{\em There are many ways to define $\<\a',\a'\>$ such that 
$V[L]$ is a vertex operator superalgebra.
For examples, one may define ${\bf h}'$ 
with $\<\a',\a'\>=1-{\overline{nk}\over n+1}$, 
where $\overline{nk}$
is the least nonnegative residue of $nk$ modulo $n+1$.
One may also define ${\bf h}'$ 
with $\<\a',\a'\>=s+{k\over n+1}$, where $s$ is any nonnegative integer.}
\er

For $\lambda\in P_{k},\;\gamma\in \C$, set
\begin{eqnarray}
W(\lambda,\gamma)=L(k,\lambda)\otimes M_{{\bf h}'}\left(1,{\gamma\over k}
\alpha'\right).
\end{eqnarray}
Since $L=\Z(h^{(1)}+\a')$,
(\ref{e4.13}) amounts to
\begin{eqnarray}
\lambda(h^{(1)})+{\gamma\over n+1}\in \Z,
\end{eqnarray}
which from (\ref{ehi}) is equivalent to
\begin{eqnarray}\label{econditionlambda}
n\lambda(\a_{1}^{\vee})+(n-1)\lambda(\a_{2}^{\vee})
+\cdots +\lambda(\a_{n}^{\vee})+\gamma\in (n+1)\Z.
\end{eqnarray}
Note that (\ref{econditionlambda}) implies $\gamma\in \Z$ 
because $\lambda\in P_{k}$.
We also see that in general, $\sigma_{W(\lambda,\gamma)}$ is of finite order 
if and only if $\gamma\in \Q$. Then Propositions \ref{p3.14} and 
\ref{pmoduleirreducibility} immediately give:

\bp{pmodule}
For $\lambda\in P_{k},\; \gamma\in \Q$, $\sigma_{W(\lambda,\gamma)}$ is of finite order 
and $W(\lambda,r)[L]$ is an irreducible $\sigma_{W(\lambda,\gamma)}$-twisted 
$A_{k}(sl(n+1))$-module. In particular, 
if (\ref{econditionlambda}) holds,
$W(\lambda,\gamma)[L]$ is an irreducible $V[L]$-module. $\;\;\;\;\Box$
\ep

{}From Theorem \ref{tregularitysemisimple} we also have:

\bp{pregularityn+1}
Let $\sigma$ be an automorphism of $A_{k}(sl(n+1))$ of finite order
which fixes $V=L(k,0)\otimes M_{{\bf h}'}(1,0)$ point-wise.
Then every $\sigma$-twisted weak $A_{k}(sl(n+1))$-module is a direct sum
of irreducible (ordinary) $\sigma$-twisted $A_{k}(sl(n+1))$-modules
$W(\lambda,\gamma)[L]$ for some $\lambda\in P_{k},\; \gamma\in \Q$ with
$\sigma =\sigma_{W(\lambda,\gamma)}$.
In particular, every weak $A_{k}(sl(n+1))$-module
is a direct sum of
$W(\lambda,\gamma)[L]$ for some $\lambda\in P_{k},\; \gamma\in \Z$
that satisfy (\ref{econditionlambda}). $\;\;\;\;\Box$
\ep

Let us consider the case $n=1$. Then we can make our results
more explicit. We have $P_{k}=\{0,1,\dots, k\}$
and $h^{(1)}={1\over 2}\a_{1}^{\vee}$.
{}From Corollary \ref{csimplecurrentfusion} we have
\begin{eqnarray}
L(k,0)^{(2mh^{(1)})}\simeq L(k,0),\;\;\;L(k,0)^{((2m+1)h^{(1)})}\simeq L(k,k)
\end{eqnarray}
for $m\in \Z$. Since
\begin{eqnarray}
V^{(m\a)}=L(k,0)^{(mh^{(1)})}\otimes M_{{\bf h}'}(1,0)^{(m\a')}
=L(k,0)^{(mh^{(1)})}\otimes M_{{\bf h}'}(1,m\alpha'),
\end{eqnarray}
it follows from Proposition \ref{pequivalencew[l]} that  
$W(i,\gamma)[L]\simeq W(i',\gamma')[L]$ if and only if
there exists $m\in \Z$ such that
\begin{eqnarray}
L(k,i')\simeq L(k,i)^{(mh^{(1)})},\;\;\;
{\gamma'\over k}={\gamma\over k}+m.
\end{eqnarray}
Recall that $[L(k,k)]\cdot [L(k,i)]=[L(k,k-i)]$.
If $m$ is even, we have $i'=i$ and $\gamma'-\gamma\in 2k\Z$.
If $m$ is odd, we have $i'=k-i$ and $\gamma'=\gamma+mk$.
Then $W(i,\gamma)[L]=W(i',\gamma')[L]$ if and only if either
$i'=i$ and $\gamma'\equiv \gamma\;\;\mod \;2k$, or
$i'=k-i$ and $\gamma'\equiv \gamma+k\;\;\mod \;2k$.
Then Propositions \ref{pmodule} and \ref{pregularityn+1}
give (cf. [FM], Proposition 3):

\bc{cfmak}
Every weak $A_{k}(sl(2))$-module is completely reducible and
\begin{eqnarray}
W(i,j)[L]\;\;\;\;\mbox{ for }0\le i\le k,\;0\le j\le k-1\;\;\;
\mbox{ with }i+j\in 2\Z
\end{eqnarray}
form a complete set of representatives of equivalence classes of irreducible
$A_{k}(sl(2))$-modules. $\;\;\;\;\Box$
\ec

Note that  $W(i,j)[L]$ was denoted by $R(i,j)$ in [FM].
Using the fusion rules for $L(k,0)$ we have the following relations
in the Verlinde algebra ${\cal{A}}(V)$:
\begin{eqnarray}
[W(i_{1},j_{1})]\cdot [W(i_{2},j_{2})]
=\sum_{i=\mbox{max}(i_{1}-i_{2},i_{2}-i_{1})}
^{\mbox{min}(i_{1}+i_{2},2k-i_{1}-i_{2})} [W(i,j_{1}+j_{2})].
\end{eqnarray}
Then in the Verlinde algebra of $A_{k}$, we have
\begin{eqnarray}
[W(i_{1},j_{1})[L]]\cdot [W(i_{2},j_{2})[L]]
=\sum_{i=\mbox{max}(i_{1}-i_{2},i_{2}-i_{1})}
^{\mbox{min}(i_{1}+i_{2},2k-i_{1}-i_{2})} [W(i,j_{1}+j_{2})[L]].
\end{eqnarray}
Note that when $j_{1}+j_{2}\ge k$, we have
$W(i,j_{1}+j_{2})[L]=W(k-i,j_{1}+j_{2}-k)[L]$.

\br{rcharacter}
{\em Set $L'=\Z\a'$. Then, as a $V$-module, 
\begin{eqnarray}
A_{k}\simeq L(k,0)\otimes V_{2L'}+L(k,k)\otimes V_{2L'+\a'}.
\end{eqnarray}
Furthermore, using more general fusion rules we get
\begin{eqnarray}
W(i,\gamma)[L]\simeq L(k,i)\otimes V_{2L'+{\gamma\over k}\alpha'}
+L(k,k-i)\otimes V_{2L'+{\gamma+k\over k}\alpha'}
\end{eqnarray}
for $i=0,\dots, k;\; \gamma\in \Q$. With this, one can easily write down
the characters of $W(i,j)[L]$ in terms of the characters of $L(k,j)$ and 
the theta functions of $2L'$.}
\er

\br{ridentification} {\em 
For $\fg=sl(n+1)$, from Corollary \ref{csimplecurrentfusion}, we have
\begin{eqnarray}
L(k,0)^{(mh^{(1)})}\simeq L(k,k\lambda_{\bar{m}})\;\;\;\mbox{ for }m\in \Z.
\end{eqnarray}
Then
\begin{eqnarray}
A_{k}(sl(n+1))\simeq
\coprod_{i=0}^{n}L(k,k\lambda_{i})\otimes V_{2(n+1)L'+i\a'}
\end{eqnarray}
as a $V$-module. More general fusion rules are needed to
express $W(\lambda,j)[L]$ explicitly.}
\er

\subsection{Generating property for the extended algebras $A_{k}(sl(n+1))$}
First, we review some properties for a general vertex operator 
superalgebra $U$. Recall Borcherds' commutator formula [B]:
\begin{eqnarray}\label{esupercommb}
[u_{m},v_{n}]_{\pm}=\sum_{i\ge 0}{m\choose i}(u_{i}v)_{m+n-i}
\end{eqnarray}
for $u,v\in U$ and $m,n\in \Z$, where $[\cdot,\cdot]_{\pm}$ refers 
to the super commutator. Thus,
the super commutator $[Y(u,z_{1}),Y(v,z_{2})]_{\pm}$ is uniquely determined by
$u_{i}v$ for $i\ge 0$. From this we have
\begin{eqnarray}\label{eweakcomm2}
(z_{1}-z_{2})^{r}[Y(u,z_{1}),Y(v,z_{2})]_{\pm}=0
\end{eqnarray}
if $r$ is a nonnegative integer such that $u_{i}v=0$ for $i\ge r$. 
For homogeneous vectors $u,v\in U$ and 
for $m\in \Z$, 
we have (cf. [FLM])
\begin{eqnarray}\label{eweightf2}
\wt (u_{m}v)=\wt u+\wt v-m-1,
\end{eqnarray}
where $\wt u$ stands for the $L(0)$-weight of $u$.

Let $U=\coprod_{n\in {1\over 2}\Z}U_{(n)}$ be 
such that $U_{(0)}=\C\;(=\C{\bf 1})$ and $U_{(n)}=0$ for $n<0$.
Then 
\begin{eqnarray}
[u_{m},v_{n}]_{+}=(u_{0}v)_{m+n-1}=\delta_{m+n,1}u_{0}v
\end{eqnarray}
for $u,v\in U_{({1\over 2})},\; m,n\in \Z$, where $u_{0}v\in U_{(0)}=\C$.
That is, the component operators $u_{m}$ for 
$u\in U_{({1\over 2})},\; m\in \Z$
give rise to a Clliford algebra.

It is well known ([B], [FLM]) that the weight-one subspace $U_{(1)}$ is a 
Lie algebra with $[u,v]=u_{0}v$ and with a symmetric invariant bilinear form 
$(u,v)=u_{1}v\in \C$. We have
\begin{eqnarray}
[u_{m},v_{n}]=(u_{0}v)_{m+n}+m\delta_{m+n,0}(u,v)
\end{eqnarray}
for $u,v\in U_{(1)},\;m,n\in \Z$. 
Then operators $u_{m}$ for $u\in U_{(1)},\; m\in \Z$
give rise to a natural representation of affine Lie algebra $\hat{U_{(1)}}$.

Now we consider $A_{k}(sl(n+1))$,
which is a vertex operator algebra
when $k$ is even and which is a vertex operator superalgebra 
when $k$ is odd. It is easy to see that
vertex operator (super)algebra $A_{k}(sl(n+1))$
is generated by $V^{(\a)}$ and $V^{(-\a)}$.
Denote by $V^{(\b)}_{low}$ the lowest $L(0)$-weight subspace 
of $V^{(\b)}$ for $\b\in L$. 
Because $V$ as a vertex operator algebra is 
generated by $\fg+\C\a'$ and both $V^{(\a)}$ and $V^{(-\a)}$ are 
irreducible $V$-modules, $A_{k}(sl(n+1))$ is furthermore generated by 
\begin{eqnarray}
S:=(\fg+\C\a')+V^{(\a)}_{low}+V^{(-\a)}_{low}.
\end{eqnarray}
Since (Corollary 2.27)
$$V^{(\a)}=L(k,k\lambda_{1})\otimes M_{{\bf h}'}(1,\a'),\;\;
V^{(-\a)}=L(k,k\lambda_{n})\otimes M_{{\bf h}'}(1,-\a'),$$
we have
\begin{eqnarray}
V^{(\a)}_{low}=L(k\lambda_{1})\otimes e^{\a'},\;\;\;
V^{(-\a)}_{low}=L(k\lambda_{n})\otimes e^{-\a'}.
\end{eqnarray}
{}From Remark \ref{rlowestweight} and (\ref{elowestweightlattice}), 
we find that 
the lowest $L(0)$-weights of $V$-modules 
$V^{(\a)}$ and $V^{(-\a)}$ are
${1\over 2}B(\a,\a)={k\over 2}$. Now we are ready to prove 
our main result of this subsection.

\bp{pgenerators}
The algebra $A_{k}(sl(n+1))$ is generated by
the subspace
$$(\fg +\C\a')+V^{(\a)}_{low}+V^{(-\a)}_{low},$$
where $V^{(\a)}_{low}=L(k\lambda_{1})\otimes e^{\a'}$ and 
$V^{(-\a)}_{low}=L(k\lambda_{n})\otimes e^{-\a'}$ are
of weight ${k\over 2}$.
Furthermore, the following relations hold:
\begin{eqnarray}
& &Y(u,z_{1})Y(v,z_{2})=(-1)^{k}Y(v,z_{2})Y(u,z_{1}),\label{esupercomm1}\\
& &Y(u',z_{1})Y(v',z_{2})=(-1)^{k}Y(v',z_{2})Y(u',z_{1})\label{esupercomm2}
\end{eqnarray}
for $u,v\in V^{(\a)}_{low},\;u',v'\in V^{(-\a)}_{low}$ and
\begin{eqnarray}
& &u_{s}v'\in V_{(k-s-1)},\label{esupercomm3}\\
& &(z_{1}-z_{2})^{k}[Y(u,z_{1}),Y(v',z_{2})]_{\pm}=0\label{esupercomm4}
\end{eqnarray}
for $u\in V^{(\a)}_{low},\; v'\in V^{(-\a)}_{low},\;s\in \Z$.
\ep

\pf First we calculate the lowest $L(0)$-weight of $V^{(m\a)}$. From 
Theorem \ref{tsimplefusionalgebra} and Corollary \ref{csimplecurrentfusion}
we have
$$V^{(m\a)}=L(k,0)^{(mh^{(1)})}\otimes M_{{\bf h}'}(1,m\a')
=L(k,k\lambda_{\bar{m}})\otimes M_{{\bf h}'}(1,m\a'),$$
where $\bar{m}$ is the least nonnegative residue of $m$ modulo $n+1$.
From Remark \ref{rlowestweight}, we see that the lowest
weight of $L(k,k\lambda_{m})$ is $\frac{\bar{m}(n+1-\bar{m})k}{2(n+1)}$. Then
the lowest weight of $V^{(m\a)}$ is
$$\frac{\bar{m}(n+1-\bar{m})k}{2(n+1)}+\frac{m^{2}k}{2(n+1)}
=\frac{\bar{m}k}{2}+\frac{(m^{2}-\bar{m}^{2})k}{2(n+1)},$$
which is at least $k$ if $|m|\ge 2$.

Let $u,v\in V^{(\a)}_{low}$. Thus $\wt u=\wt v={k\over 2}$. Then for $i\ge 0$,
$u_{i}v\in V^{(2\a)}$ and $\wt (u_{i}v)=k-i-1<k$.
Since the lowest weight of $V^{(2\b)}$ is at least $k$, we obtain
\begin{eqnarray}
u_{i}v=0\;\;\;\mbox{ for }i\ge 0.
\end{eqnarray}
Then (\ref{esupercomm1}) follows immediately from (\ref{esupercommb}).
Similarly, (\ref{esupercomm2}) holds.

(\ref{esupercomm3}) directly follows from the definition of 
the vertex operator map and the weight formula (\ref{eweightf2}).
Since $u_{s}v'\in V^{(0)}=V$ for $s\in \Z$ and 
the lowest weight of $V$ is zero,
we have
\begin{eqnarray}
u_{i}v'=0\;\;\;\mbox{ for }i\ge k.
\end{eqnarray}
Then (\ref{esupercomm4}) follows immediately from (\ref{eweakcomm2}).
$\;\;\;\;\Box$

\br{ran1} {\em 
In the case $k=1$, $L(\lambda_{1})$ is the vector representation 
of $sl(n+1)$ on $\C^{n+1}$. In this case,
the algebra $A_{1}(sl(n+1))$ is generated by  
$L(\lambda_{1})\otimes e^{\a'}+
L(\lambda_{1})^{*}\otimes e^{-\a'}$, which
generates a Clliford algebra. The algebra $A_{1}(sl(n+1))$
is exactly the spinor representation of $D_{n+1}$ [FFR], which is
isomorphic to $L(1,0)+L(1,\lambda_{1})$ as a $\hat{D}_{n+1}$-module.}
\er

\br{ran2} {\em 
When $k=2$, 
$$L(2\lambda_{1})^{*}\otimes e^{-\a'}+(\fg+\C \a')
+L(2\lambda_{1})\otimes e^{\a'}$$
is exactly the weight-one subspace of $A_{2}(sl(n+1))$
and it has a natural Lie algebra structure with the obvious $\Z$-grading. 
Using the fact that $L(2\lambda_{1})^{*}\otimes e^{-\a'}$ and 
$L(2\lambda_{1})\otimes e^{\a'}$ are non-isomorphic irreducible 
$(\fg +\C\a')$-modules one easily shows that this Lie algebra is simple
and of rank $n+1$. Consider the standard Dynkin diagram embedding of 
$sl(n+1)$ into $C_{n+1}$. Then we see
$$C_{n+1}=sl(n+1)+L(2\lambda_{1})+L(2\lambda_{n}).$$
Thus the weight one subspace of $A_{2}(sl(n+1))$ as a Lie algebra is 
isomorphic to $C_{n+1}$. (For $n=1$, this was pointed out in [FM].)
Then $A_{2}(sl(n+1))$ is a vertex operator algebra associated to
the affine Lie algebra $\hat{C}_{n+1}$.}
\er

\br{ran3} {\em For $k\ge 3$, since $\wt (u_{0}v')=k-1\ge 2$,
$[Y(u,z_{1}),Y(v',z_{2})]_{\pm}$
involves nonlinear normal ordered products of vertex operators (or fields)
$Y(a,z)$ for $a\in \fg+\C \a'$. This type of algebras are commonly referred
by physicists as nonlinear $W$-algebras.}
\er

\br{ran2'} {\em The following consideration was motivated by [GH1-2] and [Gun].
In the construction of $A_{k}$, let us define 
${\bf h}'=\C\a'$ with $\<\a',\a'\>=1+{k\over n+1}$ and keep the rest
unchanged.
Then $B(\a,\a)=1+k$. With $(L,B)$ being an integral lattice,
$V[L]$ is a vertex operator (super)algebra (cf. Remark \ref{rother}).
Furthermore,
$V[L]$ is generated by the subspace
$$V^{(-\a)}_{low}+(\fg+\C\a')+V^{(\a)}_{low}
=L(k\lambda_{1})^{*}\otimes e^{-\a'}+(\fg+\C\a')
+L(k\lambda_{1})\otimes e^{\a'},$$
where $L(k\lambda_{1})\otimes e^{\a'}$ and 
$L(k\lambda_{1})^{*}\otimes e^{-\a'}$
are of weight ${k+1\over 2}$.
In particular, when $k=2$, $V[L]$ is a vertex operator superalgebra
and $L(2\lambda_{1})\otimes e^{\a'}$ and 
$L(2\lambda_{1})^{*}\otimes e^{-\a'}$
are of weight ${3\over 2}$. In view of this and Remark \ref{ran2},
we may view $V[L]$ with $k=2$ as a superization of
the vertex operator algebra $V[L]$ with $k=2$ defined in Remark \ref{ran2}.
In [Gun], an $N=2$ vertex operator superalgebra was constructed
{}from a simple Lie algebra $\g$ equipped with a $\Z$-grading
such that $\fg_{m}=0$ for $|m|>1$. 
{}From Remark \ref{ran2}, symplectic Lie algebra $C=C_{n+1}$ 
is naturally $\Z$-graded 
with only three homogeneous subspaces being nonzero. 
A further study on the connection between $V[L]$ and the 
$N=2$ vertex operator superalgebra constructed in [Gun]
will be conducted in a future paper.}
\er



\subsection{Extended algebras $A_{k}$ of type $D_{n}$}
We consider $\fg $ of $D_{n}$ type for $n\ge 3$. 
From Corollary \ref{csimplecurrentfusion},
the equivalence classes of simple currents
$L(k,k\lambda_{1})$, $L(k,k\lambda_{n-1})$, $L(k,k\lambda_{n})$ and 
$L(k,0)$ form a group which 
is cyclic for an odd $n$ and which is isomorphic to
$\Z/2\Z\times \Z/2\Z$ for an even $n$.
We shall define the extended algebra separately for the two cases.

{}From [H] and Lemma \ref{lsimplefact} we have 
\begin{eqnarray}
& &h^{(n-1)}={1\over 2}(\a_{1}^{\vee}+2\a_{2}^{\vee}
+\cdots +(n-2)\a_{n-2}^{\vee}
+{1\over 2}n\a_{n-1}^{\vee}+{1\over 2}(n-2)\a_{n}^{\vee}),\\
& &h^{(n)}={1\over 2}(\a_{1}^{\vee}+2\a_{2}^{\vee}+\cdots +(n-2)\a_{n-2}^{\vee}
+{1\over 2}(n-2)\a_{n-1}^{\vee}+{1\over 2}n\a_{n}^{\vee})
\end{eqnarray}
and
\begin{eqnarray}
\<h^{(n-1)},h^{(n-1)}\>=\<h^{(n)},h^{(n)}\>={n\over 4}.
\end{eqnarray}
Using the relation 
$h^{(n-1)}=h^{(n)}+{1\over 2}\a_{n-1}^{\vee}-{1\over 2}\a_{n}^{\vee}$ 
we get
\begin{eqnarray}
\<h^{(n-1)},h^{(n)}\>={n-2\over 4}.
\end{eqnarray}

Case I, $n$ is odd. 

Define ${\bf h}'=\C\a'$ with $\<\a',\a'\>={3nk\over 4}$.
Set 
\begin{eqnarray}
L=\Z\a,\;\;\;\mbox{ where }\a=h^{(n)}+\a'.
\end{eqnarray}
Then
\begin{eqnarray}
L'=\Z\a',\;\;\; L''=\Z h^{(n)}.
\end{eqnarray}
Since
\begin{eqnarray}
B(\a,\a)=k\<h^{(n)},h^{(n)}\>+\<\a',\a'\>=nk,
\end{eqnarray}
$(L,B)$ is a positive definite integral lattice.
By Proposition \ref{pgeneralaffine}
(with the other assumptions being obvious)
$V[L]$ is a simple vertex operator (super)algebra.
We define $A_{k}(\fg)$ to be $V[L]$.
Then we have the following results with the same proof as that of
Propositions 4.9 and 4.10:

\bp{psonodd}
For $\fg$ of type $D_{n}$ with an odd $n$, 
the extended algebra $A_{k}(\fg)$ is regular.
Furthermore, for $\lambda\in P_{k},\; j\in \Q$, set 
$$W(\lambda,j)=L(k,\lambda)\otimes M_{{\bf h}'}(1,{j\over 3nk}\a').$$
Then any irreducible $A_{k}(\fg)$-module is isomorphic to
$W(\lambda,j)[L]$ for some $\lambda\in P_{k},\; j\in \Z$ with
\begin{eqnarray}
2\lambda(\a_{1}^{\vee})+4\lambda(\a_{2}^{\vee})+\cdots +
2(n-2)\lambda(\a_{n-2}^{\vee})
+(n-2)\lambda(\a_{n-1}^{\vee})+n\lambda(\a_{n}^{\vee})+j\in 4\Z.
\end{eqnarray}
Furthermore, $A_{k}(\fg)$ is generated by
$$(\fg+\C \a')+V_{low}^{(\a)}+V_{low}^{(-\a)},$$
where $V_{low}^{(\a)}=L(k\lambda_{n})\otimes e^{\a'}$ and
$V_{low}^{(-\a)}=L(k\lambda_{n-1})\otimes e^{-\a'}$ are of weight 
${nk\over 2}$, and the relations (\ref{esupercomm1})-(\ref{esupercomm4}) 
with $k$ being replaced by $nk$ hold.$\;\;\;\;\Box$
\ep

Case II: $n$ is even.

Define ${\bf h}'=\C\a_{1}'+\C\a_{2}'$ to be a two-dimensional vector space
with a symmetric bilinear form $\<\cdot,\cdot\>$ such that
\begin{eqnarray}
& &\<\a_{1}',\a_{1}'\>=\<\a_{1}',\a_{1}'\>={3\over 4}nk,\\
& &\<\a_{1}',\a_{2}'\>={1\over 4}k(n-2).
\end{eqnarray}
Define
\begin{eqnarray}
L=\Z\a_{1}+\Z\a_{2},
\end{eqnarray}
where
\begin{eqnarray}
\a_{1}=h^{(n-1)}+\a_{1}',\;\;\;\a_{2}=h^{(n)}+\a_{2}'.
\end{eqnarray}
Then
\begin{eqnarray}
L'=\Z\a_{1}'+\Z\a_{2}',\;\;\;L''=\Z h^{(n-1)}+\Z h^{(n)}.
\end{eqnarray}
We have
\begin{eqnarray}
& &B(\a_{1},\a_{1})=B(\a_{2},\a_{2})=kn,\\
& &B(\a_{1},\a_{2})={1\over 4}k(n-2)+{1\over 4}k(n-2)={1\over 2}k(n-2).
\end{eqnarray}
Since $n$ is even, $(L,B)$ is a positive-definite even lattice.
Clearly, $L''=\Z h^{(n-1)}+\Z h^{(n)}\subset P^{\vee}$, 
$L'$ is positive-definite, and the projection of $L$ onto $L'$ 
is one-to-one.
By Proposition \ref{pgeneralaffine}, $V[L]$ is 
a simple vertex operator algebra.
Now we define 
$A_{k}(\fg)=V[L]$, as a simple vertex operator algebra.
We just mention that this is a regular vertex operator algebra and
a set of generators and relations can be worked out similarly but with
some extra work.


\br{rdn1}
{\em Note that $L(k,k\lambda_{1})$ is a simple current of order 2 and
we have $\<h^{(1)},h^{(1)}\>=1$.
Let $V=L(k,0)$ and $L=\Z h^{(1)}$. Then in view of Corollary 3.21,
$L(k,0)+L(k,k\lambda_{1})$
has a natural simple vertex operator superalgebra structure (cf.
Remark \ref{ran1}).}
\er

\subsection{Extended vertex operator (super)algebras $A_{k}(E_{6})$}

Let $\fg$ be of type $E_{6}$. From Section 2.2, for any positive integer $k$,
$L(k,k\lambda_{1})$ and $L(k,k\lambda_{5})$ are (the only) nontrivial 
simple currents for $L(k,0)$. From [H] and Lemma \ref{lsimplefact}, we have
\begin{eqnarray}
h^{(1)}&=&{1\over 3}\left(4\a_{1}^{\vee}+3\a_{6}^{\vee}+5\a_{2}^{\vee}
+6\a_{3}^{\vee}+4\a_{4}^{\vee}+2\a_{5}^{\vee}\right),\\
h^{(5)}&=&{1\over 3}\left(2\a_{1}^{\vee}+3\a_{6}^{\vee}+4\a_{2}^{\vee}
+6\a_{3}^{\vee}+5\a_{4}^{\vee}+4\a_{5}^{\vee}\right)
\end{eqnarray}
and
\begin{eqnarray}
\<h^{(1)},h^{(1)}\>=\<h^{(5)},h^{(5)}\>={4\over 3}.
\end{eqnarray}
Define ${\bf h}'=\C\a'$ to be a one-dimensional vector space equipped with
the bilinear form $\<\cdot,\cdot\>$ such that
\begin{eqnarray}
\<\a',\a'\>={2k\over 3}.
\end{eqnarray}
Set
\begin{eqnarray}
L=\Z\a,\;\;\;\mbox{ where }\a=h^{(1)}+\a'.
\end{eqnarray}
Then 
\begin{eqnarray}
L''=\Z h^{(1)}\subset P^{\vee},\;\;\; L'=\Z\a'.
\end{eqnarray}
Since
\begin{eqnarray}
B(\a,\a)={4k\over 3}+{2k\over 3}=2k,
\end{eqnarray}
$(L,B)$ is a positive-definite even lattice. 
By Proposition \ref{pgeneralaffine} (with the other assumptions
being obvious), we have a simple 
vertex operator algebra $V[L]$.
We define $A_{k}(\fg)$ to be the vertex operator algebra $V[L]$.
For $\lambda\in P_{k},\; j\in \Q$, set
\begin{eqnarray}
W(\lambda,j)=L(k,\lambda)\otimes M_{{\bf h}'}(1,{1\over 2k}\a'),
\end{eqnarray}
an irreducible $V$-module. Then we have:

\bp{pe6}
For $\fg$ of type $E_{6}$, the extended algebra $A_{k}(\fg)$ 
is regular and any irreducible module is isomorphic to 
 $W(\lambda,j)[L]$ for some
$\lambda\in P_{k},\; j\in \Z$ with 
\begin{eqnarray}
4\lambda(\a_{1}^{\vee})+3\lambda(\a_{6}^{\vee})+5\lambda(\a_{2}^{\vee})
+6\lambda(\a_{3}^{\vee})+4\lambda(\a_{4}^{\vee})+2\lambda(\a_{5}^{\vee})
+j\in 3\Z.\;\;\;\;\Box
\end{eqnarray}
\ep

Similar to the case $\fg=sl(n+1)$,
$A_{k}(E_{6})$ as a vertex operator algebra is generated by
$$(\fg + \C \a')+V^{(\a)}_{low}+V^{(-\a)}_{low}.$$
We have
\begin{eqnarray}
& &V^{(\a)}_{low}=L(k,k\lambda_{1})_{low}\otimes e^{\a'}
=L(k\lambda_{1})\otimes e^{\a'},\\
& &V^{(-\a)}_{low}=L(k,k\lambda_{1})^{*}_{low}\otimes e^{-\a'}
=L(k\lambda_{1})^{*}\otimes e^{-\a'}.
\end{eqnarray}
Since the $L(0)$-weights of $L(k,k\lambda_{1})_{low}$ and 
$L(k,k\lambda_{1})^{*}_{low}$
are
$${1\over 2}B(h^{(1)},h^{(1)})={k\over 2}\<h^{(1)},h^{(1)}\>={2k\over 3},$$
the lowest weights of $V^{(\a)}$ and $V^{(-\a)}$ are
 ${2k\over 3}+{k\over 3}=k$.

For $m\in \Z$, the lowest $L(0)$-weight of $M_{{\bf h}'}(1,m\a')$
is ${1\over 2}\<m\a',m\a'\>={km^{2}\over 3}$. Then the lowest $L(0)$-weight
of $V^{(m\a)}$ is at least ${km^{2}\over 3}$.
If $|m|\ge 3$, the lowest $L(0)$-weight of 
$V^{(m\a)}$ is at least $3k$.
The lowest $L(0)$-weight of 
$V^{(2\a)}$ is the sum of the lowest $L(0)$-weight of $L(k,0)^{(2h^{(1)})}$
and ${4k\over 3}$. We know that $L(k,0)^{(2h^{(1)})}\simeq L(k,k\lambda_{5})$
whose lowest $L(0)$-weight is ${4k\over 3}$. Then
the lowest $L(0)$-weight of 
$V^{(2\a)}$ is ${8k\over 3}$, which is greater than $2k$.
With this information, using the same proof of
Proposition \ref{pgenerators} we immediately have:

\bp{pe6relation}
The extended vertex operator algebra $A_{k}(E_{6})$ is generated by
\begin{eqnarray}
(\fg + \C \a')+V^{(\a)}_{low}+V^{(-\a)}_{low},
\end{eqnarray}
where 
$V^{(\a)}_{low}=L(k\lambda_{1})\otimes e^{\a'}$ and
$V^{(-\a)}_{low}=L(k\lambda_{1})^{*}\otimes e^{-\a'}$
are of weight $k$.
Furthermore, the relations (\ref{esupercomm1})-(\ref{esupercomm4}) 
with $k$ being replaced by $2k$ hold.$\;\;\;\;\Box$
\ep

\br{re61}
{\em When $k=1$, 
$V^{(-\a)}_{low}+(\fg +\C\a')+V^{(\a)}_{low}$ is 
exactly the weight-one subspace of $A_{1}(\fg)$,
which is a natural Lie algebra with
\begin{eqnarray}
& &[L(\lambda_{1})\otimes e^{\a'},L(\lambda_{1})\otimes e^{\a'}]=0,\;\;\;
[L(\lambda_{1})^{*}\otimes e^{-\a'},L(\lambda_{1})^{*}\otimes e^{-\a'}]=0,\\
& &[L(\lambda_{1})\otimes e^{\a'},L(\lambda_{1})^{*}\otimes e^{-\a'}]
\subset \fg + \C \a'.
\end{eqnarray}
These relations give rise to a $\Z$-grading for the Lie algebra.
One can easily show that this Lie algebra is simple and of rank $7$.
Using the standard Dynkin diagram embedding of $E_{6}$ into $E_{7}$
we can show that it is really $E_{7}$. }
\er


\subsection{Extended vertex operator (super)algebras $A_{k}(E_{7})$}
Let $\fg$ be of type $E_{7}$. From Section 2.2, for any positive integer $k$,
$L(k,k\lambda_{6})$ is a (and the only nontrivial) simple current for $L(k,0)$.
Using [H] (Table 1 on page 69) and Lemma \ref{lsimplefact} we have
\begin{eqnarray}
h^{(6)}={1\over 2}\left(2\a_{1}^{\vee}+3\a_{7}^{\vee}+4\a_{2}^{\vee}
+6\a_{3}^{\vee}+5\a_{4}^{\vee}+4\a_{5}^{\vee}+3\a_{6}^{\vee}\right),
\;\;\;\;\;\<h^{(6)},h^{(6)}\>={3\over 2}.
\end{eqnarray}

Define ${\bf h}'=\C\a'$ to be a one-dimensional vector space equipped with
a bilinear form $\<\cdot,\cdot\>$ such that
\begin{eqnarray}
\<\a',\a'\>={k\over 2}.
\end{eqnarray}
Set 
\begin{eqnarray}
L=\Z\a,\;\;\;\mbox{ where }\a=h^{(6)}+\a'.
\end{eqnarray}
Then $L'=\Z\a'$ and $L''=\Z h^{(6)}$.
Since $B(\a,\a)={3k\over 2}+{k\over 2}=2k$, 
$(L,B)$ is a positive-definite even lattice. 
By Proposition \ref{pgeneralaffine} 
(with the other assumptions being obvious), $V[L]$ is a 
simple vertex operator algebra.
We define $A_{k}(E_{7})$ to be the simple vertex operator algebra
$V[L]$. For $\lambda\in P_{k},\; j\in \Q$, we set
\begin{eqnarray}
W(\lambda,j)=L(k,\lambda)\otimes M_{{\bf h}'}(1,{j\over k}\a').
\end{eqnarray}
In view of Theorem \ref{tregularitysemisimple} we immediately have:

\bp{pe7}
For $\fg$ of type $E_{7}$, the extended algebra 
$A_{k}(\fg)$ is regular and any irreducible module is
isomorphic to $W(\lambda,j)[L]$ for
$\lambda\in P_{k},\; j\in \Z$ with 
\begin{eqnarray}
2\lambda(\a_{1}^{\vee})+3\lambda(\a_{7}^{\vee})+4\lambda(\a_{2}^{\vee})
+6\lambda(\a_{3}^{\vee})+5\lambda(\a_{4}^{\vee})+4\lambda(\a_{5}^{\vee})+
3\lambda(\a_{6}^{\vee})+j\in 2\Z.\;\;\;\;\Box
\end{eqnarray}
\ep

The lowest weights of $V^{(\a)}$ and $V^{(-\a)}$ are ${1\over 2}B(\a,\a)=k$.
Since $[L(k,0)^{(2h^{(6)})}]=[L(k,0)]$, the lowest weights of
$V^{(2\a)}$ and $V^{(-2\a)}$ are ${1\over 2}B(2\a',2\a')=k$ also.
For $|m|\ge 3$, the lowest weight of $V^{(m\a)}$
is at least ${1\over 2}B(m\a',m\a')={m^{2}k\over 4}$, which is greater
than $2k$. With this information we immediately have:

\bp{pe7generator}
The vertex operator algebra $A_{k}(E_{7})$ is generated by
$$V_{low}^{(-2\a)}+V_{low}^{(-\a)}+(\fg+\C \a')+ V_{low}^{(\a)}+
V_{low}^{(2\a)},$$
where
\begin{eqnarray}
& &V_{low}^{(\a)}=L(k\lambda_{6})\otimes e^{\a'},\;\;\;
V_{low}^{(-\a)}=L(k\lambda_{6})\otimes e^{-\a'},\\
& &V_{low}^{(2\a)}=\C \otimes e^{2\a'},\;\;\;
V_{low}^{(-2\a)}=\C \otimes e^{-2\a'}
\end{eqnarray}
are of weight $k$. $\;\;\;\;\Box$
\ep

\br{re71}
{\em When $k=1$, 
$$\C \otimes e^{-2\a'}+L(\lambda_{6})\otimes e^{-\a'}+
(\fg+\C \a')+ L(\lambda_{6})\otimes e^{\a'}
+\C \otimes e^{2\a'}$$
is exactly the weight-one subspace of $A_{1}(\fg)$.
It is a $\Z$-graded Lie algebra with the obvious grading.
Similarly, we can show that it is $E_{8}$.}
\er


\subsection{Extended algebras of types $B_{n}$ and $C_{n}$}
For $\fg$ of type $B_{n}$, $L(k,k\lambda_{1})$ is the only
nontrivial simple current and for $\fg$ of type  $C_{n}$,
$L(k,k\lambda_{n})$ is the only nontrivial simple current. 
For $B_{n}$, from [H] and Lemma \ref{lsimplefact} we have
\begin{eqnarray}
h^{(1)}=\a_{1}^{\vee}+\cdots +a_{n-1}^{\vee}+{1\over 2}\a_{n}^{\vee},
\;\;\;\<h^{(1)},h^{(1)}\>=1
\end{eqnarray}
and for $C_{n}$ we have
\begin{eqnarray}
h^{(n)}={1\over 2}(\a_{1}^{\vee}+2\a_{2}^{\vee}+\cdots +n\a_{n}^{\vee}),
\;\;\;\;\<h^{(n)},h^{(n)}\>={n\over 2}.
\end{eqnarray}

\br{rerror}
{\em We here correct an error of [DLM2] (Examples 5.12 and 5.13)
where it was stated that $L(k,k\lambda_{n})$ was the
nontrivial simple current for $\fg$ of type $B_{n}$ and that
$L(k,k\lambda_{1})$ was the
nontrivial simple current for $\fg$ of type $C_{n}$.
For $\fg$ of type $B_{n}$,
it follows from [DLM2] or Proposition \ref{pabstract2}
with $V=L(k,0)$ and $L=\Z h^{(1)}$ that 
for any positive integral level $k$,
$L(k,0)+L(k,k\lambda_{1})$ has a natural simple 
vertex operator (super)algebra structure.
However, for $C_{n}$,
$L(k,0)+L(k,k\lambda_{n})$ is a vertex operator (super)algebra
only for a positive integral level $k$ with $nk$ being even.}
\er

For $\fg$ of type $C_{n}$, we define ${\bf h}'=\C \a'$ with 
$\<\a',\a'\>={nk\over 2}$.
Set
\begin{eqnarray}
L=\Z\a,\;\;\;\mbox{ where }\a=h^{(n)}+\a'.
\end{eqnarray}
Then $L''=\Z h^{(n)}$ and $L'=\Z \a'$. Furthermore,
\begin{eqnarray}
B(\a,\a)=k\<h^{(n)},h^{(n)}\>+\<\a',\a'\>=nk.
\end{eqnarray}
Then $(L,B)$ is a positive definite integral lattice. Hence
$V[L]$ is a simple vertex operator (super)algebra. Furthermore, 
we have:

\bp{pcn}
For $\fg$ of type $C_{n}$, $A_{k}(\fg)$ is 
regular and any irreducible $A_{k}(\fg)$-module is
isomorphic to $W(\lambda,j)[L]$ for $\lambda\in P_{k},\; j\in \Z$ with
\begin{eqnarray}
\lambda(\a_{1}^{\vee})+2\lambda(\a_{2}^{\vee})+\cdots +
n\lambda(\a_{n}^{\vee})+j\in 2\Z,
\end{eqnarray}
where
\begin{eqnarray}
W(\lambda,j)=L(k,\lambda)\otimes M_{{\bf h}'}(1,{j\over nk}\a').
\end{eqnarray}
Furthermore, $A_{k}(\fg)$ is generated by 
$$(\fg+\C \a')+V_{low}^{(\a)}+V^{(-\a)}_{low}$$
and
the the relations (\ref{esupercomm1})-(\ref{esupercomm4}) 
with $k$ being replaced by $nk$ hold, where
$V_{low}^{(\a)}$ and $V^{(-\a)}_{low}$ are of weight ${nk\over 2}$.
$\;\;\;\;\Box$
\ep

\end{document}